\documentclass[11pt]{amsart}
\usepackage[margin=1in]{geometry}
\usepackage{amsmath,amssymb,array,caption,hyperref}
\newtheorem{theorem}{Theorem}[section]
\newtheorem{lemma}[theorem]{Lemma}
\newtheorem{proposition}[theorem]{Proposition}
\newtheorem{corollary}[theorem]{Corollary}
\newtheorem{mainthm}{Theorem}

\theoremstyle{definition}
\newtheorem{definition}[theorem]{Definition}
\newtheorem{remark}[theorem]{Remark}
\newtheorem{example}[theorem]{Example}
\newtheorem{question}{Question}

\DeclareMathOperator{\Ann}{Ann}
\DeclareMathOperator{\rad}{rad}

\DeclareMathOperator{\Trd}{Trd}
\title[Prime-separated modules]{Rings in which every semiprime left ideal is an intersection of prime left ideals}
\author{Jakob Cimpri\v c}
\address{University of Ljubljana, Faculty of Mathematics and Physics, Department of Mathematics, Jadranska 21, SI-1000 Ljubljana, Slovenia, and Institute of Mathematics, Physics and Mechanics, Jadranska 19, SI-1000 Ljubljana, Slovenia}
\email{jaka.cimpric@fmf.uni-lj.si}
\thanks{The author acknowledges the financial support from the Slovenian Research and Innovation Agency (ARIS), research core funding No.\ P1-0222 and research projects No.\ J1-60011 and No.\ J1-50002.}
\date{}
\begin{document}
\begin{abstract}
We study rings in which every semiprime left ideal is an intersection
of prime left ideals, and the analogous property for submodules of
arbitrary modules. Both hold if semiprime two-sided ideals satisfy
ACC, if every principal two-sided ideal is finitely generated as
a left ideal, or if the ring is PI. Over Jacobson PI
rings, every finitely generated module is Jacobson: its prime
submodules are intersections of maximal submodules. For affine PI
algebras this extends the Amitsur--Procesi theorem to modules.

In every ring with involution, real two-sided ideals are
intersections of real prime two-sided ideals. We prove the
corresponding statement for real submodules of modules over PI rings and over
hereditary Noetherian prime rings, but show that it can
fail for real left ideals in Noetherian domains.
Replacing real prime ideals by real primitive ideals is more delicate:
even formally real affine PI algebras need not be real Jacobson, that is,
their real prime two-sided ideals need not be intersections of real
primitive ideals. We give a positivity criterion that ensures this
property. Over a real Jacobson PI ring, every finitely generated module
is real Jacobson: its real prime submodules are intersections of real
maximal submodules.

Over finitely generated free associative algebras, every finitely
generated semiprime submodule of a finite free module is a finite
intersection of finitely generated prime submodules and an intersection
of maximal submodules of finite codimension. Semiprime hulls of
finitely generated submodules are finitely generated. Evaluation at
finite-dimensional irreducible representations and nonzero vector
tuples gives a Hilbert Nullstellensatz over arbitrary fields.
For free algebras with involution, the corresponding real results are
known over $\mathbb R$ or $\mathbb C$ when the involution permutes the
generators. We extend real prime decomposition to arbitrary fields and
involutions, and give new algebraic proofs of the real Nullstellensatz
and finite generation of real radicals, without the GNS construction or
Gr\"obner bases.
\end{abstract}
\maketitle

\section{Introduction}

\subsection{Prime-separated rings}

In a commutative ring, every radical (equivalently, semiprime) ideal is
an intersection of prime ideals \cite[Tag~00E0]{Stacks}. We study the
corresponding question for one-sided ideals and, more generally,
submodules over noncommutative rings.

All rings below are associative with identity, and all modules are unital
left modules. A left ideal $I$ is \emph{semiprime} if
\[
 aRa\subseteq I\quad\Longrightarrow\quad a\in I\qquad(a\in R),
\]
and a proper left ideal $P$ is \emph{prime} if
\[
 aRb\subseteq P\quad\Longrightarrow\quad a\in P\ \text{or}\ b\in P
 \qquad(a,b\in R).
\]
These are the one-sided notions studied by Dauns
\cite[1.3]{Dauns} and Hansen \cite[Section~1]{Hansen}.
Every prime left ideal is semiprime, as is every intersection of
prime left ideals. The converse is the question motivating this paper;
see also \cite[Sections~1.3 and~3]{CS}.

\begin{question}\label{qu:leftdecomposition}
Is every semiprime left ideal of an associative unital ring an
intersection of prime left ideals?
\end{question}

A ring with the property in Question~\ref{qu:leftdecomposition} is
\emph{prime-separated}. Prime and semiprime submodules are defined
in Definition~\ref{def:prime-modules}; for left ideals these notions
agree with those above by Lemma~\ref{lem:module-characterizations}(4).
A module is \emph{prime-separated} if every semiprime submodule is
an intersection of prime submodules (Definition~\ref{def:prime-separation}).

Every module over a commutative ring is prime-separated
\cite[Theorem~A]{Cdec}. Rings in which every left ideal is idempotent
\cite[Lemma~2.1]{Hansen} and matrix rings over commutative rings
\cite[Theorem~3]{C22a} are also prime-separated. Our first result
gives two ring-theoretic sufficient conditions for every module to
be prime-separated.
An intersection over an empty family of submodules of $M$ is understood
to be $M$, and the zero module is the empty subdirect product.

\begin{mainthm}\label{thm:main}\label{thm:principal-ideal-separation}
Let $R$ be a ring satisfying at least one of the following conditions:
\begin{enumerate}
\item[(i)] $R$ satisfies ACC on semiprime two-sided ideals.
\item[(ii)] $RaR$ is finitely generated as a left ideal for every
$a\in R$.
\end{enumerate}
For every left $R$-module $M$, a submodule
$N\subseteq M$ is semiprime if and only if it is the intersection of
all prime submodules of $M$ containing it. Consequently, every
semiprime left $R$-module is a subdirect product of prime left
$R$-modules, and every semiprime left ideal of $R$ is the intersection
of the prime left ideals containing it.
\end{mainthm}

Both hypotheses concern $R$ only. The first holds for left Noetherian rings,
affine PI algebras, rings with left Krull dimension, and commutative
rings with Noetherian spectrum (Corollary~\ref{cor:ringclasses}).
The second holds, for example, for rings finite over their centre
and for tensor products over a commutative ring $k$ of left
Noetherian $k$-algebras with commutative $k$-algebras
(Corollary~\ref{cor:left-finite-ringclasses}). With no hypothesis on
the ring, locally Noetherian modules are prime-separated by
Corollary~\ref{cor:locally-noetherian}.

\subsection{Real prime-separated $\ast$-rings}

Let $R$ have an involution $a\mapsto a^*$. A left ideal $I$ is
\emph{real} if $\sum_i a_i^*a_i\in I+I^*$ implies every $a_i\in I$,
and \emph{weakly real} if $\sum_i a_i^*a_i\in I$ implies every
$a_i\in I$, for every finite family $a_i\in R$. The two notions
coincide for two-sided ideals (Subsection~\ref{subsec:realideals}) and
for maximal left ideals (Proposition~\ref{prop:reality-semisimplicity});
real modules are
defined in Definition~\ref{def:real-modules}. The ring $R$ is
\emph{real}, or formally real, if its zero ideal is real, and
\emph{semireal} if $-1$ is not a sum of hermitian squares; to avoid
confusion, algebras over the real numbers are called
$\mathbb R$-algebras. A ring, respectively a module, is \emph{real
prime-separated} if every real left ideal, respectively real
submodule, is an intersection of real prime ones
(Definition~\ref{def:real-prime-separation}). For two-sided ideals
this holds in every ring with involution.

\begin{mainthm}\label{thm:real-main}
\leavevmode
\begin{enumerate}
\item In every ring with involution, every real two-sided ideal is an
intersection of real prime two-sided ideals.
\item There is a left and right Noetherian domain with involution
containing a real left ideal that is not an intersection of real prime
left ideals.
It also contains a prime left ideal that is weakly real but not real.
\item Let $R$ be a hereditary Noetherian prime ring with involution.
Every real submodule of every left $R$-module is an intersection of
real prime submodules.
\end{enumerate}
\end{mainthm}

Part~(1) is Theorem~\ref{thm:real}. The domain in~(2) is the Jordan
plane $\mathbb R\langle x,y\rangle/(yx-xy-x^2)$ with $x^*=x$ and
$y^*=-y$ (Example~\ref{ex:jordan}); its weakly real but nonreal prime
left ideal answers \cite[Remark~5]{C22b} negatively. Part~(3) follows
from the Artinian-factor criterion in
Theorem~\ref{thm:artinian-real-factors}. It applies to noncommutative
Dedekind prime rings, to Ore extensions $D[x;\sigma,\delta]$ over
division rings with $\sigma$ an automorphism, and to their matrix
rings, but not to the Jordan plane, whose coefficient ring is
$\mathbb R[x]$ (Remark~\ref{rem:jordan-localization}).

In a Noetherian or PI ring with involution, semireality is equivalent
both to the existence of a proper real prime two-sided ideal and to
the existence of a nonzero real module.
Theorem~\ref{thm:semireal-goldie} proves these equivalences whenever
$R/P$ is left Goldie for every involution-invariant prime two-sided
ideal $P$. Semireality need not imply the existence of a real simple
module, or equivalently a real maximal left ideal
(Proposition~\ref{prop:real-simple}).

\subsection{Jacobson properties and PI rings}

An affine commutative algebra over a field $k$ is \emph{Jacobson}: every
prime ideal is an intersection of maximal ideals, and these have finite
$k$-codimension \cite[Tag~00FV]{Stacks}. If $k$ is formally real and the
involution is trivial, every real prime ideal is an intersection of real
maximal ideals \cite[Proposition~6.6 and Remark~6.7]{Cnull}; here and
below, \emph{real maximal} means real and maximal.

A two-sided ideal is \emph{primitive} if it is the annihilator of a
simple module, and \emph{real primitive} if it is the annihilator of a
real simple module. In commutative rings these are the maximal and the
real maximal ideals. More generally, in PI
rings primitive ideals are precisely maximal two-sided ideals by
Kaplansky's theorem \cite[Part~B, Theorem~4.13]{DF}, and real primitive
ideals are precisely real maximal two-sided ideals
(Proposition~\ref{prop:PI-real-primitive}).

A ring $R$ is \emph{Jacobson} if every prime two-sided ideal is an
intersection of primitive ideals, and \emph{left Jacobson} if every
prime left ideal is an intersection of maximal left ideals. A left
$R$-module is \emph{Jacobson} if every prime submodule is an intersection
of maximal submodules.
For rings with involution, the corresponding real versions are
obtained by adding the word ``real'' throughout these definitions.
For arbitrary rings, left Jacobson implies Jacobson, and real left
Jacobson implies real Jacobson. The converses can fail, but they hold
for PI rings.

\begin{mainthm}\label{thm:loc-main}
Let $R$ be a PI ring.
\begin{enumerate}
\item Every left $R$-module is prime-separated and, if $R$ has an
involution, real prime-separated.
\item $R$ is Jacobson if and only if it is left Jacobson, if and only if
every finitely generated left $R$-module is Jacobson.
\item If $R$ has an involution, it is real Jacobson if and only if it is
real left Jacobson, if and only if every finitely generated left
$R$-module is real Jacobson.
\item If $R$ is an affine algebra over a field $k$, then in every
finitely generated left $R$-module, every semiprime submodule is an
intersection of maximal submodules, and every maximal submodule has
finite $k$-codimension.
\item Suppose that $R$ is an Azumaya algebra whose centre $C$ is affine
over a formally real field $k$, and that the involution is $k$-linear.
If the reduced trace $\Trd_{R/C}(a^*a)\,1_R$ is a sum of hermitian
squares in $R$ for every $a\in R$, then $R$ is real Jacobson.
\end{enumerate}
\end{mainthm}

Part~(4) is the module form of the Amitsur--Procesi theorem
\cite[Part~B, Theorem~11.7]{DF}, which says that affine PI algebras
are Jacobson. Part~(5) extends to Azumaya algebras the
Procesi--Schacher sufficient trace criterion for central simple algebras
\cite[p.~404]{PS}. Real Jacobsonness can fail for affine PI algebras:
Example~\ref{ex:twistedreal} is a formally real affine PI
$\mathbb R$-algebra with no real primitive ideals at all. For affine
PI $k$-algebras with a $k$-linear involution,
Theorem~\ref{thm:realPI}(ii) gives a more general criterion, which is harder
to verify: the involution on the central simple ring of quotients of every real
prime factor must be positive at some ordering of its fixed centre.
This criterion is sufficient but not necessary for $R$ to be real
Jacobson (Remark~\ref{rem:positivity-not-necessary}). It also ensures
that every real primitive ideal is the kernel of a unital
$*$-homomorphism $R\to(M_n(K),{}^{\mathsf T})$ with $K$ a real closure
of $k$ (Corollary~\ref{cor:positive-real-primitive}).

Part~(1) is proved by central localization (Theorem~\ref{thm:PI}).
Parts~(2) and~(3) follow from Theorem~\ref{thm:jacobson-transfer},
which applies to rings whose finitely generated prime modules embed in
products of their annihilator quotients; PI rings have this property
(Lemma~\ref{lem:PIembedding}).
Parts~(4) and~(5) are Corollaries~\ref{cor:affine-PI-maximal}
and~\ref{cor:azumaya-trace-square}, respectively.

\subsection{Free algebras}

For free associative algebras, we must restrict ourselves to finitely
generated submodules of finite free modules. We obtain a version of Hilbert's Nullstellensatz:

\begin{mainthm}\label{thm:free-main}
Let $F=k\langle x_1,\ldots,x_g\rangle$ over an arbitrary field $k$, with
$g<\infty$.
\begin{enumerate}
\item Every finitely presented semiprime left $F$-module is a finite
direct sum of prime modules. Consequently, a finitely generated
submodule $N\subseteq F^n$ is semiprime if and only if it is a finite
intersection of finitely generated prime submodules.
\item A finitely generated submodule $N\subseteq F^n$ is semiprime if and
only if it is an intersection of finitely generated maximal submodules.
Equivalently, every finitely presented semiprime left $F$-module is a
subdirect product of finite-dimensional simple modules. A maximal
submodule of $F^n$ is finitely generated exactly when it has finite
codimension over $k$.
\item If $N$ is a finitely generated submodule of $F^n$, its smallest
semiprime hull $\rad N$ is finitely generated and $\rad N/N$ is
finite-dimensional. Generators of $\rad N$ can be computed by a
terminating algorithm when $k=\mathbb Q$ or $k$ is a finite field.
\end{enumerate}
\end{mainthm}

The finitely generated maximal submodules in part~(2) can be described
explicitly. Choose a finite-dimensional irreducible representation
$\pi:F\to\operatorname{End}_k(V)$ and vectors $v_1,\ldots,v_n\in V$,
not all zero. Define the $F$-linear map
\begin{equation}\label{eq:intro-evaluation}
 \varepsilon_{\pi,v}:F^n\longrightarrow V,\qquad
 (f_1,\ldots,f_n)\longmapsto\sum_j\pi(f_j)v_j.
\end{equation}
The kernels $L_{\pi,v}=\ker\varepsilon_{\pi,v}$ are exactly the
maximal submodules of finite $k$-codimension
(Proposition~\ref{prop:free-maximal}).

It is classical that every finitely generated left ideal in a free
associative algebra is an intersection of left ideals of finite
codimension: see Rosenmann--Rosset \cite[Theorem~4.3]{RR} (they work with right ideals), and Bergman's
Nullstellensatz in
\cite{HM}. Those statements allow representations that need
not be irreducible. Restricting to irreducible representations gives the
semiprime hull instead of the original ideal, a question raised in
\cite[Section~2.3]{CS}. Theorem~\ref{thm:free-main} answers this question
and treats submodules of $F^n$ at the same time. For $n=1$, its maximal-submodule assertion gives a version of the Jacobson
property for finitely generated prime left ideals, without asserting
that the free algebra is left Jacobson.

\subsection{Free algebras with involution}

Parts~(2) and~(3) of Theorem~\ref{thm:free-real-main} follow from
\cite{CHMN,HKN,Nelson}, as outlined in
Remark~\ref{rem:earlier-real-nullstellensatze}. Our contribution here
is that we extend part~(1) to arbitrary fields and involutions and
that we reprove part~(2) and the finite-generation and degree-bound
assertions of part~(3) by algebraic methods, without using the GNS
construction or Gr\"obner bases.

\begin{mainthm}\label{thm:free-real-main}
Let $k$ be a field with involution and let
$F=k\langle x_1,\ldots,x_g\rangle$, with $g<\infty$, carry an involution
extending that of $k$.
\begin{enumerate}
\item Every finitely presented real left $F$-module is a finite
direct sum of real prime modules. Consequently, every finitely
generated real submodule of $F^n$ is a finite intersection of finitely
generated real prime submodules.
\end{enumerate}
For~(2) and~(3), suppose that $k=\mathbb R$ with the identity
involution or $k=\mathbb C$ with complex conjugation, and that $*$
permutes $x_1,\ldots,x_g$.
\begin{enumerate}
\setcounter{enumi}{1}
\item \textup{\cite{CHMN,HKN,Nelson}}
A finitely generated submodule $N\subseteq F^n$ is real if and only
if it is an intersection of finitely generated real maximal submodules.
Equivalently, every finitely presented real left $F$-module is a
subdirect product of finite-dimensional real simple modules. A real
maximal submodule of $F^n$ is finitely generated exactly when it has
finite codimension over $k$. These submodules are kernels of the maps
in \eqref{eq:intro-evaluation} with $\pi$ an irreducible $*$-representation.
\item \textup{\cite{CHMN,HKN,Nelson}}
The real radical $\sqrt[\mathrm{re}]{N}$ of every finitely generated
submodule $N\subseteq F^n$ is finitely generated. If $N$ is generated
in degrees at most $d\ge1$, its real radical is generated by the given
generators and at most $n\sum_{j=0}^{d-1}g^j$ further vectors of degree
less than $d$. It can be computed by a terminating algorithm for real
algebraic coefficients, or complex coefficients with real algebraic
real and imaginary parts.
\end{enumerate}
\end{mainthm}

\subsection{Related work and organization}

Applications of this theory to polynomial and
group rings will be developed in separate work.

Sections~\ref{sec:main-proof} and~\ref{sec:involution} prove
Theorems~\ref{thm:main} and~\ref{thm:real-main}, respectively.
Section~\ref{sec:localization} proves Theorem~\ref{thm:loc-main}.
Sections~\ref{sec:free} and~\ref{sec:freereal} prove
Theorems~\ref{thm:free-main} and~\ref{thm:free-real-main}, except for
the algorithmic assertions. Appendix~\ref{sec:algorithms} gives the
prime-radical algorithm and recalls the established real-radical
algorithm.

\section{Prime separation}\label{sec:main-proof}

Theorem~\ref{thm:main} gives two conditions on the ring under which
every semiprime submodule is an intersection of prime submodules.
We also prove that Noetherian modules over arbitrary rings are
prime-separated.

\subsection{Prime and semiprime submodules}\label{subsec:prime-semiprime-submodules}

Prime and semiprime modules extend the corresponding notions for
left ideals. The definitions below let us treat all modules in a
single statement and recover the left-ideal result as a special case.

\begin{definition}[Prime and semiprime modules]\label{def:prime-modules}
For a left $R$-module $L$, its \emph{annihilator} is
\[
 \Ann_R(L)=\{a\in R:a\ell=0\text{ for every }\ell\in L\}.
\]
It is a two-sided ideal of $R$.
A nonzero left $R$-module $Q$ is \emph{prime} if every nonzero submodule
$L\subseteq Q$ has $\Ann_R(L)=\Ann_R(Q)$. A left $R$-module $Q$ is
\emph{semiprime} if
\[
 L\cap\Ann_R(L)Q=0\qquad\text{for every submodule }L\subseteq Q.
\]
A submodule $N\subseteq M$ is \emph{prime} or \emph{semiprime} if
$M/N$ is a prime or semiprime module, respectively.
\end{definition}

Lemma~\ref{lem:module-characterizations} gives practical tests for these definitions and
checks that they recover the corresponding notions for left ideals.

\begin{lemma}[Characterizations]\label{lem:module-characterizations}
Let $Q,M$ be left $R$-modules and $N\subseteq M$ a submodule.
\begin{enumerate}
\item A nonzero $Q$ is prime if and only if
\[
 aRq=0\quad\Longrightarrow\quad q=0\ \text{or}\ aQ=0
 \qquad(a\in R,\ q\in Q).
\]
\item The module $Q$ is semiprime if and only if
\[
 q\in\Ann_R(Rq)Q\quad\Longrightarrow\quad q=0\qquad(q\in Q),
\]
or, equivalently, $L\subseteq\Ann_R(L)Q$ implies $L=0$ for every
submodule $L\subseteq Q$.
\end{enumerate}
In particular:
\begin{enumerate}
\setcounter{enumi}{2}
\item The submodule $N$ is prime if and only if $N\ne M$ and
\[
 aRm\subseteq N\quad\Longrightarrow\quad m\in N\ \text{or}\ aM\subseteq N
 \qquad(a\in R,\ m\in M).
\]
It is semiprime if and only if, for every finite sum
$m=\sum_i a_i m_i$ with $a_i\in R$ and $m_i\in M$,
\[
 a_iRm\subseteq N\ \text{for all }i\quad\Longrightarrow\quad m\in N.
\]
\item A left ideal is prime (respectively semiprime) as a submodule
of ${}_R R$ if and only if it is prime (respectively semiprime) as a
left ideal.
\end{enumerate}
\end{lemma}

\begin{proof}
(1) The criterion says that every nonzero cyclic submodule has
annihilator $\Ann_R(Q)$. This is equivalent to the same assertion
for all nonzero submodules.

(2) Semiprimeness gives both implications. Conversely, the condition
on submodules, applied to $H=L\cap\Ann_R(L)Q$, gives $H=0$ because
$H\subseteq\Ann_R(L)Q\subseteq\Ann_R(H)Q$.
The cyclic condition also gives $H=0$:
every $q\in H$ belongs to $\Ann_R(Rq)Q$, since
$\Ann_R(L)\subseteq\Ann_R(Rq)$.

(3) Apply (1) and (2) to $M/N$. For the semiprime criterion, expand
$q\in\Ann_R(Rq)(M/N)$ as $q=\sum_i a_i(m_i+N)$ with
$a_iRq=0$, and take $m=\sum_i a_i m_i$.

(4) The prime criterion in (3) implies the left-ideal criterion.
Conversely, if $aRb\subseteq I$ and $b\notin I$, applying the
left-ideal criterion to $ac,b$ for every $c\in R$ gives $aR\subseteq I$.
For semiprimeness, take $m=a\cdot1$ in (3) to obtain the left-ideal
criterion. Conversely, if $m=\sum_i a_i b_i$ and every $a_iRm\subseteq I$,
then $mRm\subseteq\sum_i a_iRm\subseteq I$, so that criterion gives
$m\in I$.
\end{proof}

For modules, the term \emph{semiprime} has several meanings in the
literature. We use the same notion as in \cite{Cdec}
(Lemma~\ref{lem:module-characterizations}(2)); Jir\'asko calls it
\emph{pseudo-semiprimeness} \cite[Proposition~1.2(vi)]{Jirasko}.
See also the discussion preceding \cite[Proposition~5.5]{Lomp}.
By Lemma~\ref{lem:module-characterizations}(2), our notion implies
semiprimeness in the sense of Dauns \cite{DaunsModules}, namely
$aRam=0\Rightarrow am=0$. The converse fails even over commutative
rings \cite[p.~93]{Lomp}.

The following closure properties let us build semiprime modules
from simpler ones.

\begin{lemma}[Prime and semiprime modules]\label{lem:basic-prime-modules}
Let $R$ be a ring.
\begin{enumerate}
\item Every simple $R$-module is prime. Every prime $R$-module is semiprime.
\item Every nonzero submodule of a prime module is a prime module.
\item Every submodule of a semiprime module is a semiprime module.
\item Every subdirect product of semiprime modules is semiprime.
In particular, products and direct sums of semiprime modules, and hence
of prime modules, are semiprime.
\end{enumerate}
\end{lemma}

\begin{proof}
(1) A simple module has no nonzero submodule other than itself, so
it is prime. If $Q$ is prime and $0\ne L\subseteq Q$, then
$\Ann_R(L)Q=\Ann_R(Q)Q=0$; for $L=0$ the semiprimeness condition
is automatic. Thus $Q$ is semiprime.

(2) If $0\ne L\subseteq Q$ with $Q$ prime, then every nonzero
$H\subseteq L$ satisfies $\Ann_R(H)=\Ann_R(Q)=\Ann_R(L)$,
so $L$ is prime.

(3) If $N\subseteq Q$ with $Q$ semiprime, then for every $L\subseteq N$,
\[
 L\cap\Ann_R(L)N\subseteq L\cap\Ann_R(L)Q=0,
\]
so $N$ is semiprime.

(4) Let $V$ be a subdirect product of semiprime modules $V_\lambda$,
with quotient maps $f_\lambda:V\twoheadrightarrow V_\lambda$.
For any submodule $L\subseteq V$, we have
$\Ann_R(L)\subseteq\Ann_R(f_\lambda(L))$, hence
\[
 f_\lambda\bigl(L\cap\Ann_R(L)V\bigr)
 \subseteq f_\lambda(L)\cap\Ann_R(f_\lambda(L))V_\lambda=0.
\]
Since the maps jointly detect zero, $L\cap\Ann_R(L)V=0$.
Products and direct sums are subdirect products of their factors;
the assertion for prime modules also uses (1).
\end{proof}

Primeness and semiprimeness of submodules are defined by the
corresponding properties of their factor modules. Therefore
Lemma~\ref{lem:basic-prime-modules} immediately implies
Lemma~\ref{lem:basic-prime-submodules}.

\begin{lemma}[Prime and semiprime submodules]\label{lem:basic-prime-submodules}
Let $M$ be an $R$-module.
\begin{enumerate}
\item Every maximal submodule of $M$ is prime. Every prime submodule
of $M$ is semiprime.
\item Every intersection of semiprime submodules of $M$ is semiprime;
in particular, every intersection of prime submodules is semiprime.
\item If $N\subseteq L\subseteq M$, then $L$ is a prime submodule
(respectively a semiprime submodule) of $M$ if and only if $L/N$ is a
prime submodule (respectively a semiprime submodule) of $M/N$.
\end{enumerate}
\end{lemma}

\begin{proof}
For (1), apply Lemma~\ref{lem:basic-prime-modules}(1) to quotient
modules. For (2), $M/\bigcap_\lambda N_\lambda$ is a subdirect
product of the semiprime modules $M/N_\lambda$, so
Lemma~\ref{lem:basic-prime-modules}(4) applies; the assertion for
prime submodules also uses (1). For (3), use
$(M/N)/(L/N)\cong M/L$.\qedhere
\end{proof}

Annihilators connect prime and semiprime modules with the familiar
two-sided ideal notions.

\begin{lemma}[Prime and semiprime annihilators]\label{lem:prime-semiprime-annihilators}
\leavevmode\par
\begin{enumerate}
\item The annihilator of a prime module is a prime two-sided ideal.
\item The annihilator of a semiprime module is a semiprime two-sided
ideal. In particular, the annihilator of every submodule of a semiprime
module is semiprime.
\end{enumerate}
\end{lemma}

\begin{proof}
(1) Let $Q$ be prime and $T=\Ann_R(Q)$. This ideal is proper since
$Q\ne0$. If two-sided ideals $A,B$ satisfy $AB\subseteq T$ and
$B\not\subseteq T$, then $BQ\ne0$. Primeness gives
$\Ann_R(BQ)=T$, and $A(BQ)=0$ therefore implies $A\subseteq T$.
Thus $T$ is prime.

(2) Let $Q$ be semiprime and let $A$ be a two-sided ideal with
$A^2Q=0$. Put $H=AQ$. Then $AH=0$, so
$H=AQ\subseteq\Ann_R(H)Q$.
Lemma~\ref{lem:module-characterizations}(2) gives $H=0$,
hence $A\subseteq\Ann_R(Q)$. Thus $\Ann_R(Q)$ is semiprime.
The assertion for submodules follows from
Lemma~\ref{lem:basic-prime-modules}(3).
\end{proof}

The separation argument will use the following way of extending a
prime module.

\begin{lemma}[Essential extensions of prime modules]\label{lem:essential-prime}
Let $U$ be an essential submodule of $M$, meaning that
$U\cap H\ne0$ for every nonzero submodule $H\subseteq M$.
If $U$ is prime and $\Ann_R(M)=\Ann_R(U)$, then $M$ is prime.
In particular, every essential extension of a faithful prime module
is faithful and prime.
\end{lemma}

\begin{proof}
For $0\ne H\subseteq M$, essentiality gives $H\cap U\ne0$.
Primeness of $U$ therefore yields
\[
 \Ann_R(M)\subseteq\Ann_R(H)
 \subseteq\Ann_R(H\cap U)=\Ann_R(U)=\Ann_R(M).
\]
Thus every nonzero submodule of $M$ has the same annihilator as $M$,
so $M$ is prime. If $U$ is faithful, then
$\Ann_R(M)\subseteq\Ann_R(U)=0$, proving the last assertion.
\end{proof}

\subsection{Prime separation}\label{subsec:prime-separation}

We now formulate prime separation for modules. The observations and
the Separation Lemma below apply over any ring.

\begin{definition}[Prime separation]\label{def:prime-separation}
A left $R$-module $M$ is \emph{prime-separated} if every semiprime
submodule is an intersection of prime submodules. Equivalently,
for every semiprime $N\subseteq M$ and $m\notin N$, some prime
submodule contains $N$ and misses $m$.
A ring is \emph{prime-separated} if its regular left module is.
We state explicitly when the property holds for all left modules.
\end{definition}

When every $R$-module is prime-separated, a semiprime module can
be recovered from its prime quotients: it embeds naturally as a
subdirect product of them. Taking $M={}_R R$ gives the left-ideal
assertion of Theorem~\ref{thm:main}
(Lemma~\ref{lem:module-characterizations}(4)). The smallest semiprime
left ideal containing a left ideal $J$, called its \emph{semiprime
hull}, is then the intersection of all prime left ideals containing
$J$. Conversely, intersections of prime submodules are always
semiprime (Lemma~\ref{lem:basic-prime-submodules}(2)).

Write $\beta(R)$ for the \emph{lower nilradical} (or two-sided
prime radical) of $R$, the intersection of all prime two-sided
ideals. It is also the intersection of all prime left ideals:
by Lemma~\ref{lem:prime-semiprime-annihilators}(1), every prime
left ideal $P$ contains the prime two-sided ideal $\Ann_R(R/P)$,
while every prime two-sided ideal is prime as a left ideal.

Part~(3) of the next lemma uses Theorem~\ref{thm:main}, proved
below independently of this lemma.

\begin{lemma}[Basic constructions]\label{lem:prime-separation-constructions}
\leavevmode
\begin{enumerate}
\item Prime-separated rings are closed under quotients and finite
products. If $I\subseteq\beta(R)$ is a two-sided ideal, then $R$
is prime-separated if and only if $R/I$ is. The same assertions hold
for the property that every left module is prime-separated.
\item For any rings $A,B$ and $(A,B)$-bimodule $M$, each of these
two properties holds for
$T=\left(\begin{smallmatrix}A&M\\0&B\end{smallmatrix}\right)$
if and only if it holds for both $A$ and $B$.
Triangular constructions preserve hypothesis~(i) of
Theorem~\ref{thm:main}, but need not preserve~(ii).
\item If $R/\beta(R)$ satisfies~(ii), every left $R$-module is
prime-separated. This condition is preserved by quotients, finite
products, matrix rings and triangular constructions.
\end{enumerate}
\end{lemma}
\begin{proof}
(1) Restriction along a quotient map preserves primeness and
semiprimeness, giving closure under quotients. If
$I\subseteq\beta(R)$ and $L\subseteq V$ is semiprime, then
$IV\subseteq L$, since $\Ann_R(V/L)$ is semiprime by
Lemma~\ref{lem:prime-semiprime-annihilators}(2).
Thus prime and semiprime submodules correspond under
$V\to V/IV$, proving both equivalences. For products, modules
decompose componentwise and prime quotients come from one factor.

(2) The off-diagonal ideal of $T$ is square-zero and its quotient
is $A\times B$, so (1) applies. Semiprime two-sided ideals of $T$
likewise correspond to pairs of such ideals of $A,B$, proving
preservation of~(i). Example~\ref{ex:biregular-principal-ideals}
shows that~(ii) need not be preserved.

(3) The first assertion follows from Theorem~\ref{thm:main}
and (1). To prove closure, first suppose that $R$ itself
satisfies~(ii). Quotients inherit finite left generating sets, and
finite products satisfy~(ii) componentwise.
For $S=M_n(R)$ and $b=(b_{ij})\in S$,
$SbS=M_n(I)$, where $I=\sum_{i,j}Rb_{ij}R$ is finitely generated
as a left ideal of $R$. If $I=\sum_{\ell=1}^d Rt_\ell$, then
\[
 M_n(I)=\sum_{\ell=1}^d\sum_{j=1}^n S(t_\ell E_{1j}),
\]
so matrix rings also preserve~(ii). Apply these observations to
$R/\beta(R)$, using the identities
$\beta(M_n(R))=M_n(\beta(R))$ and
$T/\beta(T)\cong A/\beta(A)\times B/\beta(B)$,
and the fact that every semiprime factor of $R$ is a quotient of
$R/\beta(R)$.
\end{proof}

Thus Question~\ref{qu:leftdecomposition} reduces to semiprime rings.

\begin{lemma}[Separation]\label{lem:separate}
Let $N$ be a semiprime submodule of $M$ and $x\in M\setminus N$.
Suppose that $(N+Rx)/N$ contains a nonzero prime submodule. Then
some prime submodule of $M$ contains $N$ but not $x$.
\end{lemma}

\begin{proof}
Choose $N\subsetneq L\subseteq N+Rx$ with $L/N$ prime and put
$T=\Ann_R(L/N)$. Semiprimeness of $M/N$ gives
$(L/N)\cap T(M/N)=0$, or equivalently $L\cap(N+TM)=N$.

Choose $P$ maximal among submodules containing $N+TM$ and satisfying
$P\cap L=N$. Such a $P$ exists by Zorn's lemma: $N+TM$ satisfies
these conditions, and unions of chains preserve them. If $x\in P$,
then $L\subseteq N+Rx\subseteq P$, contradicting $P\cap L=N$ and
$L\ne N$. Thus $x\notin P$.

The image $\overline L=(L+P)/P\cong L/N$ is essential in $M/P$:
for $P\subsetneq U\subseteq M$, maximality gives
$U\cap L\supsetneq N$, so $U/P$ meets $\overline L$ nontrivially.
Moreover, $\overline L$ is prime with annihilator $T$, and
$TM\subseteq P$ gives
\[
 T\subseteq\Ann_R(M/P)\subseteq\Ann_R(\overline L)=T.
\]
Lemma~\ref{lem:essential-prime} makes $M/P$ prime, as required.
\end{proof}

\subsection{ACC on two-sided semiprime ideals}\label{subsec:ACC-semiprime-ideals}

An ascending chain condition on semiprime two-sided ideals is
enough to ensure prime separation for every module. This covers,
in particular, left Noetherian rings and affine PI algebras. We also
show that locally Noetherian modules are prime-separated over any ring.

\begin{lemma}[Existence of prime submodules]
\label{lem:prime-submodule-exists}
Let $Q$ be semiprime and $0\ne q\in Q$. The cyclic submodule $Rq$
contains a nonzero prime submodule in each of the following cases:
\begin{enumerate}
\item $R$ satisfies ACC on semiprime two-sided ideals;
\item $Rq$ contains a nonzero uniform submodule, in particular if
$Rq$ is Noetherian. Here a nonzero module is \emph{uniform} if any
two nonzero submodules intersect nontrivially.
\end{enumerate}
\end{lemma}
\begin{proof}
(1) Each $\Ann_R(L)$, $0\ne L\subseteq Rq$, is semiprime by
Lemma~\ref{lem:prime-semiprime-annihilators}(2). Choose $L$ with
$\Ann_R(L)$ maximal. For $0\ne H\subseteq L$, maximality gives
$\Ann_R(H)=\Ann_R(L)$, so $L$ is prime.

(2) Let $0\ne U\subseteq Rq$ be uniform. If $Rq$ is Noetherian,
such a submodule exists: otherwise successively choosing two
disjoint nonzero submodules and continuing inside one of them
produces an infinite direct sum whose partial sums contradict ACC.
For $0\ne H\subseteq U$,
\[
 H\cap\Ann_R(H)U\subseteq H\cap\Ann_R(H)Q=0.
\]
Since $H\ne0$ and both $H$ and $\Ann_R(H)U$ are submodules of
the uniform module $U$, their zero intersection forces
$\Ann_R(H)U=0$. Hence $\Ann_R(H)\subseteq\Ann_R(U)$; the reverse
inclusion follows from $H\subseteq U$. Thus every nonzero submodule
of $U$ has the same annihilator as $U$, so $U$ is prime.
\end{proof}

\begin{proof}[Proof of Theorem~\ref{thm:main} under hypothesis~(i)]
Let $N\subseteq M$ be semiprime and $x\notin N$.
Lemma~\ref{lem:prime-submodule-exists}(1), applied to $Q=M/N$ and
$q=x+N$, supplies the hypothesis of Lemma~\ref{lem:separate}.
Thus a prime submodule contains $N$ and misses $x$; conclude by
the remarks after Definition~\ref{def:prime-separation}.
\end{proof}

\begin{corollary}[Locally Noetherian modules]
\label{cor:locally-noetherian}
Over an arbitrary ring, every locally Noetherian module, that is,
every module whose cyclic submodules are Noetherian, is
prime-separated. This applies in particular to Noetherian modules
and to finite-dimensional modules over any algebra over a field.
\end{corollary}
\begin{proof}
Let $N$ be a semiprime submodule of such a module $M$ and
$x\notin N$. The cyclic module $(N+Rx)/N$, being a quotient of
$Rx$, is Noetherian. Apply Lemma~\ref{lem:prime-submodule-exists}(2)
to $M/N$ and then Lemma~\ref{lem:separate}.
\end{proof}

\begin{corollary}[Classes of rings]\label{cor:ringclasses}
Each of the following rings satisfies ACC on semiprime two-sided
ideals and hence the conclusions of Theorem~\ref{thm:main} for every
left module:
\begin{enumerate}
\item a left Noetherian ring;
\item an affine PI algebra over a field;
\item a ring with left Krull dimension;
\item a commutative ring whose spectrum is Noetherian.
\end{enumerate}
\end{corollary}

\begin{proof}
(1) Every ascending chain of two-sided ideals is an ascending chain
of left ideals, so it stabilizes.

(2) Affine PI algebras satisfy ACC on semiprime two-sided ideals by
\cite[Corollary~1, p.~182]{Chiba}.

(3) Apply \cite[Theorem~7.6]{GR} to $R^{\mathrm{op}}$, which has
right Krull dimension when $R$ has left Krull dimension. The two rings
have the same two-sided ideals and the same semiprime two-sided ideals,
so the resulting ACC holds for $R$.

(4) In a commutative ring, semiprime ideals are radical ideals.
The correspondence $I\mapsto V(I)$ reverses inclusion and bijects
radical ideals with closed subsets of $\operatorname{Spec}R$.
An ascending chain of radical ideals therefore gives a descending
chain of closed subsets, which stabilizes by Noetherianity of the
spectrum. Hence the original chain stabilizes as well.

Theorem~\ref{thm:main} now applies in each case.
\end{proof}

Theorem~\ref{thm:PI} will extend prime separation to all PI rings.

\subsection{Principal two-sided ideals as left ideals}\label{subsec:principal-ideal-separation}

The role of hypothesis~(ii) of Theorem~\ref{thm:main} is to ensure
that unions of nonempty chains of semiprime submodules are semiprime.
We first show
that this chain-union property implies prime separation, and then
verify it under~(ii). The hypothesis covers every commutative ring
and tensor products over a commutative ring $k$ of left Noetherian
$k$-algebras with arbitrary commutative $k$-algebras.

\begin{lemma}[Annihilator quotients]\label{lem:annihilator-quotients}
Let $Q$ be a semiprime module and $I$ a two-sided ideal. Then
$Q/(0:_Q I)$ is semiprime, where $(0:_Q I)=\{q\in Q:Iq=0\}$.
\end{lemma}
\begin{proof}
Put $K=(0:_Q I)$. For a submodule $L\supseteq K$, let
$J=\Ann_R(L/K)$, so $JL\subseteq K$. Then
$IJ\subseteq\Ann_R(L)$, which is semiprime by
Lemma~\ref{lem:prime-semiprime-annihilators}(2). Since
$(JI)^2\subseteq J(IJ)I$, we also have $JI\subseteq\Ann_R(L)$.
Thus $JIL=0$, so $J\subseteq\Ann_R(IL)$.
For $x\in L\cap(JQ+K)$, write $x=y+k$ with $y\in JQ$ and
$k\in K$. Then $Ix=Iy\subseteq IJQ\subseteq JQ$, while
$Ix\subseteq IL$. Hence
\[
 Ix\subseteq IL\cap JQ
 \subseteq IL\cap\Ann_R(IL)Q=0.
\]
Thus $L\cap(JQ+K)\subseteq K$. Since
$J(Q/K)=(JQ+K)/K$, this is the semiprimeness condition for $Q/K$.
\end{proof}

The preceding quotient construction lets us show that a semiprime
submodule maximal subject to avoiding an element is prime.

\begin{lemma}[Maximal semiprime submodules]
\label{lem:maximal-semiprime-submodules}
Let $x\in M$, and let $N$ be maximal among the semiprime
submodules of $M$ avoiding $x$. Then $N$ is prime.
Consequently, $M$ is prime-separated whenever unions of nonempty
chains of semiprime submodules of $M$ are semiprime.
\end{lemma}
\begin{proof}
Put $Q=M/N$ and $q=x+N$. For $0\ne H\subseteq Rq$, set
$J=\Ann_R(H)$. The submodule $(0:_Q J)$ contains $H\ne0$ and
is semiprime in $Q$ by Lemma~\ref{lem:annihilator-quotients}.
Its preimage in $M$ is a semiprime submodule strictly larger than
$N$, so maximality forces it to contain $x$. Thus $Jq=0$.
Since $J$ is two-sided, $JRq\subseteq Jq=0$, giving
$\Ann_R(H)\subseteq\Ann_R(Rq)$. The reverse inclusion is
automatic; hence $Rq$ is prime. Lemma~\ref{lem:separate} supplies
a prime submodule containing $N$ and avoiding $x$, which must
equal $N$ by maximality.

Under the chain-union hypothesis, Zorn's lemma gives such an $N$
containing any prescribed semiprime submodule and avoiding any
element outside it. The first assertion then proves prime separation.
\end{proof}

\begin{proof}[Proof of Theorem~\ref{thm:main} under hypothesis~(ii)]
By Lemma~\ref{lem:maximal-semiprime-submodules}, it suffices to
show that a nonempty chain $(N_\lambda)$ of semiprime submodules
has semiprime union $N$. Suppose
\[
 x-\sum_{i=1}^r a_i m_i\in N,
 \qquad a_iRx\subseteq N\quad(1\le i\le r).
\]
Since the elements $a_is$, with $s\in R$, generate $Ra_iR$ as a left
ideal, hypothesis~(ii) gives finite generating sets
$Ra_iR=\sum_{j=1}^{d_i}Ra_i s_{ij}$.
All the finitely many elements
$x-\sum_i a_i m_i$ and $a_i s_{ij}x$ lie in one $N_\lambda$.
Then $Ra_iRx=\sum_j Ra_i s_{ij}x\subseteq N_\lambda$, hence
$a_iRx\subseteq N_\lambda$. Thus in $M/N_\lambda$ the image
of $x$ belongs to
$\Ann_R(R(x+N_\lambda))(M/N_\lambda)$. Semiprimeness gives
$x\in N_\lambda\subseteq N$.
The elementwise criterion proves that $N$ is semiprime.
\end{proof}

\begin{corollary}[Classes of rings]
\label{cor:left-finite-ringclasses}
Each of the following rings $R$ satisfies hypothesis~(ii) of
Theorem~\ref{thm:main}:
\begin{enumerate}
\item a commutative ring;
\item a ring finitely generated as a module over its centre
(e.g.\ an Azumaya algebra);
\item a ring generated by a left Noetherian subring and a
central subring (e.g.\ $B\otimes_k S$, with $k$ a commutative ring,
$B$ a left Noetherian $k$-algebra and $S$ a commutative $k$-algebra);
\item a biregular ring (Example~\ref{ex:biregular-principal-ideals}).
\end{enumerate}
Consequently, every left module over these rings is
prime-separated.
\end{corollary}
\begin{proof}
For (1), $RaR=Ra$. For (2), if $b_1,\ldots,b_d$ generate $R$
over $Z(R)$, then
$RaR=\sum_{j=1}^d Rab_j$.
For (3), let $B$ be the left Noetherian subring and $C$ the central
subring. Centrality gives $R=BC$, the additive span of products
$bc$. Write $a=\sum_{i=1}^d b_i c_i$. Then
$BaB\subseteq\sum_i Bc_i$, a Noetherian left $B$-module.
Hence $BaB=\sum_{j=1}^e Bu_j$ for finitely many $u_j$, and
$RaR=C(BaB)=\sum_{j=1}^e Ru_j$.
In the tensor example, use the images of $B\otimes1$ and
$1\otimes S$; the former is left Noetherian as a quotient of $B$.
For (4), $RaR=Re$ for a central idempotent $e$.
\end{proof}

The two hypotheses are independent: $k[t_1,t_2,\ldots]$ satisfies
\textup{(ii)} but not~\textup{(i)}, whereas the PI ring of
Remark~\ref{rem:jacobson-centre} satisfies~\textup{(i)} but
not~\textup{(ii)}.

\begin{example}[FC-group crossed products]\label{ex:FC-crossed-products}
Let $B$ be a left Noetherian ring and $G$ a group whose conjugacy
classes are finite (an \emph{FC-group}). Every module over a crossed
product $T=B*G$ is prime-separated. Indeed, write
$T=\bigoplus_{g\in G}Bu_g$ with homogeneous units $u_g$.
For $a\in T$ with finite support $F$, the set
$F^G=\bigcup_{g\in G}gFg^{-1}$ is finite, and the $B$-subbimodule
\[
 V=\sum_{g\in G}B(u_g^{-1}au_g)B
 \subseteq\bigoplus_{h\in F^G}Bu_h
\]
is finitely generated on the left. Since $TaT=TV$,
Theorem~\ref{thm:main} applies under hypothesis~(ii). The equality
follows by moving each right factor $u_g$ to the left:
$abu_g=u_g(u_g^{-1}au_g)(u_g^{-1}bu_g)$ for $b\in B$.
This includes skew Laurent extensions in possibly infinitely many
commuting variables associated with commuting automorphisms of $B$;
no finite-order assumption on the automorphisms is needed.
Further prime-separation and Jacobson properties of group rings
will be studied in a separate paper.
\end{example}

\begin{example}[Biregular rings]\label{ex:biregular-principal-ideals}
A ring is \emph{biregular} if every principal two-sided ideal is
generated by a central idempotent. Let $k$ be a field. For instance, in
\[
 A=k\mathbf1+\bigoplus_{n\ge2}M_n(k)
 \subseteq\prod_{n\ge2}M_n(k)
\]
every $a$ has finite or cofinite support, and $AaA=Ae_a$ for the
corresponding central idempotent $e_a$. The ring $A$ is not PI, as it
maps onto every $M_n(k)$, and fails~(i), as the finite sums of its
matrix components form a strictly increasing chain of semiprime
ideals.

Now put $U=\left(\begin{smallmatrix}k&A\\0&A\end{smallmatrix}\right)$.
Its quotient $U/\beta(U)\cong k\times A$ satisfies~(ii), so every
$U$-module is prime-separated by
Lemma~\ref{lem:prime-separation-constructions}. Yet $U$ is not PI
and fails~(i), since $A$ is a quotient. It also fails~(ii):
$Ue_{11}U\cong k\oplus A$ is infinite-dimensional over $k$, with
left $U$-action through the upper-left scalar entry.
\end{example}

\section{Rings with involution}\label{sec:involution}

The first three subsections prove the three parts of
Theorem~\ref{thm:real-main}: real two-sided prime separation, the
Jordan-plane counterexample, and real prime separation over hereditary
Noetherian prime rings. The last subsection studies semireality
and the existence of real maximal left ideals. Section~\ref{sec:localization}
then develops the PI case.

\subsection{Real two-sided ideals}\label{subsec:realideals}

Let $R$ be a ring with involution $a\mapsto a^*$. A left ideal $I$ is
\emph{real} if $\sum_i a_i^*a_i\in I+I^*$ implies every $a_i\in I$,
and \emph{weakly real} if $\sum_i a_i^*a_i\in I$ implies every
$a_i\in I$, for every finite family $a_i\in R$. The ring $R$ is
\emph{real} if its zero ideal is real.
Since $I\subseteq I+I^*$, reality implies weak reality.
For two-sided ideals the two notions agree: if $I$ is weakly real
and $a\in I$, then $(a^*)^*a^*=aa^*\in I$ gives $a^*\in I$, so
$I=I^*$ and $I+I^*=I$. They also agree for maximal left ideals
(Proposition~\ref{prop:reality-semisimplicity}), but not in general
(Examples~\ref{ex:jordan} and~\ref{ex:freeweak}).
Weak reality is the notion called ``real'' in
\cite[\S1.3]{C22b}; whether the two definitions agree for all left
ideals was asked in \cite[Remark~5]{C22b}.

In a real $*$-ring,
the annihilator $T=\Ann_R(J)$ of any left ideal
$J$ is a real two-sided ideal. Indeed, if $\sum_i a_i^*a_i\in T$, then
for every $j\in J$,
\[
 \sum_i(a_ij)^*(a_ij)=j^*\Bigl(\sum_i a_i^*a_i\Bigr)j=0.
\]
Reality of $R$ gives $a_iJ=0$, so every $a_i\in T$. Thus $T$ is
weakly real and, being two-sided, real.

\begin{theorem}[Real two-sided ideals]\label{thm:real}
In every $*$-ring, every real two-sided ideal is an intersection of real
prime two-sided ideals.
\end{theorem}
\begin{proof}
For two-sided ideals, reality and weak reality agree; hence a union
of a chain of real two-sided ideals is real. Let $I$ be a real
two-sided ideal and $a\notin I$. By
Zorn's lemma there is a real two-sided ideal $P\supseteq I$ maximal with
respect to $a\notin P$. Put $B=R/P$ and $b=a+P\ne0$. The ring $B$ is
real, and by maximality every nonzero real two-sided ideal of $B$
contains $b$. A real ring is semiprime: if $cBc=0$, then
$(c^*c)^*(c^*c)=c^*(cc^*c)=0$, so $c=0$.

Suppose that $UV=0$ for nonzero two-sided ideals $U,V$ of $B$. Then
$(VU)^2\subseteq V(UV)U=0$, so $VU=0$. Put $T=\Ann_B(U)$ and
$S=\Ann_B(T)$. Since $B$ is real, the preceding annihilator
observation shows that $T$ and $S$ are real self-adjoint two-sided
ideals. Now $V\subseteq T$, and $TU=0$ gives
$(UT)^2\subseteq U(TU)T=0$, so $UT=0$ and $U\subseteq S$. Hence $S$ and
$T$ are nonzero and both contain $b$. Since $S=S^*$, we get
$b^*b\in ST=0$, so $b=0$, a contradiction. Thus $P$ is a real prime
two-sided ideal containing $I$ and not $a$.
\end{proof}

\begin{remark}[Related results]\label{rem:nilpotent-enveloping}
For Noetherian complex $*$-algebras, Schmitt and Sch\"otz
\cite[Theorem~28]{SS} prove the finite-intersection version of
Theorem~\ref{thm:real}. For a finite-dimensional nilpotent real Lie
algebra $\mathfrak g$, put
$U^*(\mathfrak g)=U(\mathfrak g)\otimes_{\mathbb R}\mathbb C$,
with $x^*=-x$ for $x\in\mathfrak g$. Their stronger result
\cite[Theorem~7 and Remark~33]{SS} expresses every real two-sided
ideal of $U^*(\mathfrak g)$ as an intersection of kernels of
surjective $*$-homomorphisms onto Weyl algebras. These kernels are
real maximal and completely prime two-sided ideals, since Weyl
algebras are formally real simple domains.
\end{remark}

\subsection{The Jordan plane}\label{subsec:jordan}

Theorem~\ref{thm:real} does not extend to left ideals, even in a
Noetherian domain. The following example also distinguishes weak
reality from reality, answering \cite[Remark~5]{C22b} negatively.

\begin{example}[The Jordan plane]\label{ex:jordan}
Let
\[
 A=\mathbb R\langle x,y\rangle/(yx-xy-x^2),\qquad
 x^*=x,\quad y^*=-y,
\]
and put $b=y+x^2$, $P=Ab$ and $I=Axb$. Then $A$ is a left and right
Noetherian domain, and:
\begin{enumerate}
\item $P$ is prime and weakly real but not real;
\item $I$ is real, but every real prime left ideal containing $I$
contains $x$, while $x\notin I$. Hence $I$ is
not an intersection of real prime left ideals.
\end{enumerate}
\end{example}

\begin{proof}
We use a polynomial model and an integral pairing to establish
the reality assertions, and normality of $x$ to control the real
prime quotients.

\emph{Structure.} The Ore extension $A=\mathbb R[x][y;\delta]$, with
$\delta=x^2\,d/dx$, is a left and right Noetherian domain, and $*$
preserves the relation. Since $yx=x(y+x)$, the element $x$ is normal,
so $(x)=Ax=xA$. The quotient $A/(x)\cong\mathbb R[y]$, with $y^*=-y$,
is a real domain: evaluating $\sum_i\bar a_i(-y)\bar a_i(y)=0$ at
imaginary points gives $\sum_i|\bar a_i(is)|^2=0$. Hence $(x)$ is a
real prime left ideal.

\emph{A model of $A/P$.} The operators
\[
 Xp=tp,\qquad Yp=t^2p'-t^2p
\]
on $\mathbb R[t]$ satisfy $[Y,X]=X^2$ and $b\cdot1=0$. Since $b$ is
monic in $y$, every element of $A$ is congruent modulo $P$ to a
polynomial in $x$, so $a+P\mapsto a\cdot1$ identifies $A/P$ with
$\mathbb R[t]$. In particular $x\notin P$. Moreover $I=P\cap(x)$:
clearly $I\subseteq P\cap xA$, and if $ub\in(x)$, then reduction
modulo $(x)$ gives $\bar u\,y=0$ in $\mathbb R[y]$, so $u=u'x$ and
$ub=u'xb\in I$. Thus $x\notin I$.

\emph{The pairing used to prove reality.} Let $\nu(f)$ denote the order of vanishing
of $f\in\mathbb R[t]$ at $0$, with $\nu(0)=\infty$; both $X$ and $Y$
raise it by at least one.
Whenever the integral converges, put
\[
 \langle f,g\rangle=\int_0^\infty f(t)g(t)\,t^{-2}e^{-2t}\,dt .
\]
If $\nu(f)+\nu(g)\ge1$ and $w$ is a nonempty word in $X,Y$, then both
sides of
\begin{equation}\label{eq:jordanadj}
 \langle w^*f,g\rangle=\langle f,wg\rangle
\end{equation}
converge and are equal. For $X$ this is clear. For $Y$, since
$t^2\cdot t^{-2}e^{-2t}=e^{-2t}$, integration by parts gives
$\langle Yf,g\rangle+\langle f,Yg\rangle=[fge^{-2t}]_0^\infty=0$.
Induction on the length of $w$ proves~\eqref{eq:jordanadj}, because each
generator raises $\nu$. By linearity,~\eqref{eq:jordanadj} also holds
with $w$ replaced by any element of $A$ with zero constant term.

\emph{Weak reality and failure of reality.} Since
$b+b^*=2x^2=x^*x+x^*x\in P+P^*$ and $x\notin P$, the
left ideal $P$ is not real. Suppose $\sum_ia_i^*a_i\in P$ and put
$p_i=a_i\cdot1$; then $\sum_ia_i^*\cdot p_i=0$.
The augmentation $x,y\mapsto0$ kills $P=Ab$, so
$\sum_ic_i^2=0$, where $c_i$ is the constant term of $a_i$.
Thus every $c_i=0$, each $a_i$ is a combination of nonempty words
and $\nu(p_i)\ge1$, and~\eqref{eq:jordanadj} with $(f,g)=(p_i,1)$ gives
\[
 0=\sum_i\langle a_i^*p_i,1\rangle=\sum_i\int_0^\infty
 p_i(t)^2t^{-2}e^{-2t}\,dt.
\]
Hence every $p_i=0$, that is, $a_i\in P$. So $P$ is weakly real.

\emph{Primeness.} Every nonzero two-sided ideal
$H$ of $A$ contains a power of $x$. A nonzero element of $H$ of least
$y$-degree has degree zero, since otherwise its commutator with $x$
has smaller $y$-degree and nonzero leading coefficient. Thus
$H\cap\mathbb R[x]=(f)$ is nonzero. This contraction is stable under
$\delta=x^2d/dx$, so $f$ divides $x^2f'$. In characteristic zero this
forces $f$ to be a nonzero scalar times a power of $x$.
Now $x$ acts injectively on $A/P\cong\mathbb R[t]$. The annihilator
of any nonzero submodule therefore cannot contain a power of $x$,
and hence is zero. Thus $A/P$, and therefore $P$, is prime.

\emph{(2), reality.} Suppose $S=\sum_ia_i^*a_i\in I+I^*$. Since
$I+I^*\subseteq(x)$ and $(x)$ is real, every $a_i\in(x)$; hence
$\nu(p_i)\ge1$ for $p_i=a_i\cdot1$, and $a_i$ has no constant term.
Write $S=uxb+b^*xv^*$ and $h=xv^*\cdot1$, so $\nu(h)\ge1$ and
$S\cdot1=b^*\cdot h$ because $b\cdot1=0$. Applying~\eqref{eq:jordanadj}
twice,
\[
 \sum_i\int_0^\infty p_i^2\,t^{-2}e^{-2t}\,dt
 =\sum_i\langle a_i^*p_i,1\rangle=\langle b^*h,1\rangle
 =\langle h,b\cdot1\rangle=0 .
\]
Thus every $p_i=0$, so $a_i\in P\cap(x)=I$, and $I$ is real.

\emph{(2), separation.} Let $J\supseteq I$ be a real prime left ideal,
$H=A/J$ and $q=1+J\ne0$. By normality, $K=\ker(x\colon H\to H)$ is a
submodule. If $K=0$, then $xbq=0$ gives $bq=0$, so $b\in J$ and
$x^*x+x^*x=b+b^*\in J+J^*$. Reality gives $x\in J$, so $xq=0$ and
$q\in K=0$, a contradiction. Hence $K\ne0$, and primeness gives
$x\in\Ann_A(K)=\Ann_A(H)$. Therefore $(x)\subseteq J$.
\end{proof}

\subsection{Hereditary Noetherian prime rings}\label{subsec:HNP-real-prime}\label{subsec:realbasic}

The Jordan plane of Example~\ref{ex:jordan} is the Ore extension
$A=\mathbb R[x][y;x^2d/dx]$, so real prime separation can fail
with one Ore variable over a commutative polynomial ring. For an
Ore extension over a division ring, every nonzero proper two-sided
ideal has an Artinian quotient. The same quotient property holds
for hereditary Noetherian prime rings, including noncommutative
Dedekind prime rings. For real ideals these quotients are semisimple,
which supplies the splitting argument below. In the Jordan plane,
the real ideal $(x)$ instead has the non-Artinian quotient
$A/(x)\cong\mathbb R[y]$.

We prove the separation statement for arbitrary modules, beginning
with the definitions and basic properties needed for the argument.

For a left module $Q$, give its conjugate additive group $\overline Q$
the right action $\bar q a=\overline{a^*q}$. The tensor product
$\overline Q\otimes_R Q$ is an abelian group; its balancing relation is
\[
 (\bar q a)\otimes v=\bar q\otimes(av)\qquad(a\in R).
\]

\begin{definition}[Real and weakly real modules]\label{def:real-modules}
An $R$-module $Q$ is \emph{real} if, for every finite family $q_i\in Q$,
\[
 \sum_i\bar q_i\otimes q_i=0\quad\text{in }\overline Q\otimes_R Q
 \quad\Longrightarrow\quad q_i=0\text{ for all }i,
\]
and \emph{weakly real} if, for every $q\in Q$ and every finite family
$a_i\in R$,
\begin{equation}\label{eq:weakreal}
 \Bigl(\sum_i a_i^*a_i\Bigr)q=0\quad\Longrightarrow\quad
 a_iq=0\text{ for all }i.
\end{equation}
A submodule $N\subseteq M$ is \emph{real} or \emph{weakly real} if
$M/N$ is real or weakly real, respectively, and \emph{real prime} if
$M/N$ is real and prime.
\end{definition}

Every real module is weakly real, by the identity
\[
 \sum_i\overline{a_iq}\otimes a_iq
 =\bar q\otimes\Bigl(\sum_i a_i^*a_i\Bigr)q .
\]

The next proposition checks that these module definitions recover
the definitions for left ideals.

\begin{proposition}[Cyclic real modules]\label{prop:cyclicreal}
Let $I$ be a left ideal of $R$.
\begin{enumerate}
\item $I$ is real if and only if the left module $R/I$ is real.
\item $I$ is weakly real if and only if the left module $R/I$ is weakly
real. In that case $I$ is semiprime as a left ideal and $R/I$ is
semiprime as a module.
\end{enumerate}
\end{proposition}
\begin{proof}
(1) The map $\bar a\otimes b\mapsto a^*b$ induces an additive
isomorphism $\overline{R/I}\otimes_R R/I\cong R/(I^*+I)$, where the
right-hand side is a quotient of additive groups, with inverse
$c+(I^*+I)\mapsto\overline{1+I}\otimes(c+I)$.
The inverse vanishes on $I$, and also on $I^*$, since for $i\in I$,
\[
 \overline{1+I}\otimes(i^*+I)
 =\overline{i+I}\otimes(1+I)=0.
\]
The balancing relation shows that the two maps are inverse.
The isomorphism sends $\sum_i\bar a_i\otimes a_i$ to
$\sum_ia_i^*a_i+(I+I^*)$.

(2) If $I$ is weakly real and $(\sum_i a_i^*a_i)b\in I$, left
multiplication by $b^*$ gives $\sum_i(a_ib)^*(a_ib)\in I$, hence every
$a_ib\in I$; the converse follows by taking $b=1$. If $aRa\subseteq I$,
then $(a^*a)^*(a^*a)=a^*(aa^*a)\in I$, and weak reality gives
successively $a^*a\in I$ and $a\in I$. Thus $I$ is semiprime, and so is
$R/I$ by Lemma~\ref{lem:module-characterizations}(4).
\end{proof}

The following closure and annihilator properties will be used in
the separation arguments.

\begin{lemma}[Basic properties]\label{lem:apprealbasic}\label{lem:weakann}
Real and weakly real modules are closed under submodules, products and
restriction of scalars along $*$-homomorphisms. Real modules are
semiprime and weakly real. In a
weakly real module, the annihilator of every submodule is a real
self-adjoint two-sided ideal.
\end{lemma}
\begin{proof}
For real modules, the closure assertions follow by mapping tensor
relations into the ambient tensor product, into each coordinate, or along
the ring map. For weakly real modules they follow directly from
\eqref{eq:weakreal}.
Weak reality of real modules was noted above. Let $Q$ be weakly real and
$T=\Ann_R(L)$ for a submodule $L$. If $a\in T$, then $a^*L\subseteq L$ gives $aa^*L=0$, so weak reality
gives $a^*L=0$; thus $T=T^*$, and applying~\eqref{eq:weakreal} to each
element of $L$ proves reality of $T$. For real $Q$ and
$z=\sum_jt_jq_j\in L\cap TQ$, we get
$\bar z\otimes z=\sum_j\overline{t_j^*z}\otimes q_j=0$, so $z=0$;
this is semiprimeness of $Q$.
\end{proof}

\begin{remark}[Unions of real submodules]\label{rem:real-routes}
Unions of nonempty chains of real submodules are real. Indeed, if
$N=\bigcup_\lambda N_\lambda$, a vanishing hermitian sum
$\sum_i\bar m_i\otimes m_i$ in
$\overline{M/N}\otimes_R(M/N)$ already vanishes for some
$M/N_\lambda$, since tensor products commute with directed colimits.
Reality of $M/N_\lambda$ then forces every $m_i\in N_\lambda\subseteq N$.
Thus the chain-union step in the maximality argument of
Lemma~\ref{lem:maximal-semiprime-submodules} is automatic for
real submodules. Maximal real submodules avoiding a prescribed
element exist, but they need not be prime, as the cyclic Noetherian
real module $A/I$ of Example~\ref{ex:jordan} shows.
\end{remark}

Intersections of real submodules are real, since the quotient embeds
in the product of the real quotients. Their intersection defines the
\emph{real radical} $\sqrt[\mathrm{re}]{N}$ of a submodule $N$.

The proof uses only the Artinian quotient property for nonzero
real ideals, as the following criterion shows.

\begin{theorem}[Real prime separation from Artinian factors]
\label{thm:artinian-real-factors}
Let $R$ be a ring with involution such that $R/T$ is left Artinian
for every nonzero real two-sided ideal $T$. Then every real
submodule of every left $R$-module is an intersection of real prime
submodules.
\end{theorem}
\begin{proof}
It suffices to separate each nonzero element $q$ of a real module $M$.
By Remark~\ref{rem:real-routes} and Zorn's lemma, choose a real
submodule $P$ maximal subject to $q\notin P$. Put $Q=M/P$.
If $Q=Q_1\oplus Q_2$ with both summands nonzero, then each $Q_i$
is a real submodule of $Q$, since $Q/Q_i$ is isomorphic to a
submodule of the real module $Q$. At least one $Q_i$ misses $q+P$,
and its preimage contradicts maximality of $P$. Thus $Q$ is
indecomposable.

If $Q$ is not prime, choose $0\ne H\subsetneq Q$ with
$T=\Ann_R(H)\supsetneq\Ann_R(Q)$; in particular, $T\ne0$.
By Lemma~\ref{lem:apprealbasic}, $T$ is real and self-adjoint.
By hypothesis $R/T$ is left Artinian. Being real, it is semiprime
and therefore semisimple Artinian.
Semiprimeness of $Q$ gives $H\cap TQ=0$, so $H$ embeds in $Q/TQ$.
Both are modules over the semisimple ring $R/T$, so this embedding
has a retraction. Composing with $Q\to Q/TQ$ makes $H$ a nonzero
proper direct summand of $Q$, a contradiction. Hence $Q$ is prime,
and $P$ is a real prime submodule missing $q$.
\end{proof}

Every quotient of a hereditary Noetherian prime ring by a nonzero
two-sided ideal is Artinian \cite[Lemma~2(1)]{TuganbaevHNP}.
The theorem therefore proves Theorem~\ref{thm:real-main}(3), including
the case of noncommutative Dedekind prime rings. In particular,
$D[x;\sigma,\delta]$, with $D$ a division ring, $\sigma$ an
automorphism and $\delta$ a $\sigma$-derivation, is a left and right
principal ideal domain, hence a hereditary Noetherian prime ring.
Here $xa=\sigma(a)x+\delta(a)$ for $a\in D$.
Its matrix rings are hereditary
Noetherian prime as well. Thus the theorem applies to all these
rings with any involution.

\begin{remark}[Localizing the Jordan plane]\label{rem:jordan-localization}
Passing to division coefficients does not preserve the real left
ideal in Example~\ref{ex:jordan}. Indeed, the derivation and
involution extend to
\[
 \widetilde A=\mathbb R(x)[y;x^2d/dx],\qquad x^*=x,\quad y^*=-y,
\]
so Theorem~\ref{thm:artinian-real-factors} applies to $\widetilde A$.
But $x$ is now invertible, and the real left ideal $I=Axb$ extends
to $\widetilde A I=\widetilde A b$, where $b=y+x^2$.
This left ideal is proper by degree in $y$, but it is not real:
\[
 b+b^*=2x^2=x^*x+x^*x
 \in\widetilde A b+(\widetilde A b)^*.
\]
Reality would force the unit $x$ into $\widetilde A b$.
Thus the theorem over division coefficients cannot be applied to
this extended ideal to obtain real prime separation of $I$ in $A$.
\end{remark}

Hereditary Noetherian prime rings need not be real Jacobson, even
when they are Jacobson.

\begin{example}[The integers]\label{ex:hnp-not-real-jacobson}
With the identity involution, $\mathbb Z$ is real and, as a principal
ideal domain, hereditary Noetherian prime. It is Jacobson: all nonzero
prime ideals are maximal, and $(0)=\bigcap_{p\text{ prime}}p\mathbb Z$.
However, no maximal ideal $p\mathbb Z$ is real, since the sum of $p$
copies of $1^*1$ lies in it while $1\notin p\mathbb Z$. Thus $(0)$ is
real prime, but there are no real primitive ideals, so $\mathbb Z$
is not real Jacobson.
\end{example}

\subsection{Semireality and real maximal left ideals}\label{subsec:semireal}

Write $\Sigma R^*R$ for the set of finite sums of hermitian squares
$a^*a$. A $*$-ring $R$ is \emph{semireal} if $-1\notin\Sigma R^*R$.
Theorem~\ref{thm:real} separates elements from real two-sided ideals;
the following result gives a criterion for such an ideal to exist.

Recall that an element of a ring $B$ is \emph{regular} if it is
neither a left nor a right zero divisor. When the set $S$ of all regular
elements satisfies the left Ore condition, the localization $S^{-1}B$
is called the \emph{classical left quotient ring}, or \emph{classical
left localization}, of $B$; its elements have the form $s^{-1}b$ with
$s\in S$ and $b\in B$. Goldie's theorem guarantees its existence for
a semiprime left Goldie ring and says that it is semisimple Artinian;
if $B$ is prime, it is simple Artinian \cite[2.3.6]{MR}.
Every semisimple Artinian ring $Q$ is von Neumann regular: for every
$a\in Q$ there exists $b\in Q$ such that $aba=a$.

\begin{theorem}[Semireality]\label{thm:semireal-goldie}
Let $R$ be a nonzero $*$-ring such that $R/J$ is left Goldie for every
proper involution-invariant prime two-sided ideal $J$ (for example, a
left or right Noetherian $*$-ring by Goldie's theorem
\cite[Chapter~2]{MR}, or a PI $*$-ring by Posner's theorem
\cite[Part~B, Section~6]{DF}).
The following conditions are equivalent:
\begin{enumerate}
\item $R$ is semireal.
\item $R$ has a proper real prime two-sided ideal.
\item $R$ has a proper real two-sided ideal.
\item $R$ has a proper real left ideal.
\item There exists a nonzero real left $R$-module.
\end{enumerate}
These conditions need not imply the existence of a real maximal
left ideal, even when $R$ is a formally real affine PI algebra.
\end{theorem}

\begin{proof}
We first prove (1)$\Rightarrow$(2). A \emph{quadratic module} on a
$*$-ring $B$ is a subset of $\operatorname{Sym}(B)=\{b:b^*=b\}$ containing
$1$, closed under addition and under $m\mapsto x^*mx$ for $x\in B$;
it is proper if it does not contain $-1$. By Zorn's lemma,
$\Sigma R^*R$ is contained in a maximal proper quadratic module $M$.
By \cite[Theorem~1]{Cmax}, there is a proper involution-invariant prime
two-sided ideal $J$ with
\[
 J\cap\operatorname{Sym}(R)=M\cap(-M).
\]
Its image $\bar M$ in $B=R/J$ is proper and has
$\bar M\cap(-\bar M)=0$. Indeed, if $m,n\in M$ have opposite images,
then $m+n\in M\cap(-M)$, so $-m=n-(m+n)\in M$ and $m\in J$.
If $-1$ were the image of $m\in M$, then
$1+m\in J\cap\operatorname{Sym}(R)=M\cap(-M)$ would give
$-1=m-(1+m)\in M$.

By hypothesis $B$ is prime and left Goldie; its involution makes it
right Goldie as well. Goldie's theorem gives a simple Artinian classical
quotient ring $Q$. The left and right classical quotient rings coincide,
so the involution extends to $Q$. By
\cite[Proposition~2]{Cmax},
\[
 \{q\in\operatorname{Sym}(Q):c^*qc\in\bar M
       \text{ for some regular }c\in B\}
\]
is a proper quadratic module on $Q$. Choose a maximal proper quadratic
module $P$ containing it. Again by \cite[Theorem~1]{Cmax},
$P\cap(-P)$ is the symmetric part of a proper two-sided ideal of $Q$.
Simplicity gives $P\cap(-P)=0$.

We show that $Q$ is real. If $\sum_i a_i^*a_i=0$, each $a_i^*a_i$
lies in $P\cap(-P)$, hence is zero. For $a=a_i$, von Neumann regularity
gives $b\in Q$ with $aba=a$. Put $p=ab$ and $e=p+p^*-pp^*$.
Then $p^2=p$, $p^*p=0$, $a=pa$, and direct expansion gives
$e^*=e$, $e^2=e$ and $ep=p$. Consequently
\[
 e=e^*e\in P,\qquad -e=(p-p^*)^*(p-p^*)\in P.
\]
Thus $e=0$, then $p=ep=0$ and $a=pa=0$. Therefore $Q$, and its
$*$-subring $B$, are real, so $J$ proves~(2).

The implications (2)$\Rightarrow$(3)$\Rightarrow$(4) are immediate.
If $I$ is a proper real left ideal and $-1=\sum_i a_i^*a_i$, then
$1^*1+\sum_i a_i^*a_i=0\in I$ forces $1\in I$, a contradiction.
This proves (4)$\Rightarrow$(1).
The implication (4)$\Rightarrow$(5) follows from
Proposition~\ref{prop:cyclicreal}. Conversely, the annihilator of a
nonzero real module is a proper real two-sided ideal by
Lemma~\ref{lem:apprealbasic}, proving (5)$\Rightarrow$(3).
The final assertion follows from Example~\ref{ex:twistedreal}.
\end{proof}

Only (1)$\Rightarrow$(2) uses the Goldie hypothesis.

To study real maximal left ideals, we first establish the relevant
facts about semisimple modules and weak reality.

\begin{proposition}[Reality and semisimplicity]\label{prop:reality-semisimplicity}
Let $R$ be a unital $*$-ring.
\begin{enumerate}
\item If $R$ is real and semisimple Artinian, every left $R$-module
is real.
\item A semisimple left $R$-module is real if and only if it is
weakly real.
\end{enumerate}
Consequently, a maximal left ideal, or a finite intersection of
maximal left ideals, is real if and only if it is weakly real.
This equivalence can fail for infinite intersections
(Example~\ref{ex:freeweak}). Here maximal means maximal among all
proper left ideals.
\end{proposition}
\begin{proof}
(1) The regular module of $R$ is real by
Proposition~\ref{prop:cyclicreal}. For a simple module $U$, choose
$0\ne u\in U$. The map $\pi:R\to U$, $a\mapsto au$, is surjective.
Semisimplicity gives $R=\ker\pi\oplus L$, so $U\cong L\subseteq R$
is real by Lemma~\ref{lem:apprealbasic}. Every module is a direct sum
of such real simple modules, hence is real by the same lemma, since
the direct sum embeds in their product.

(2) Only the implication from weak reality to reality needs proof.
First let $L$ be a weakly real maximal left ideal, and suppose
$\sum_i a_i^*a_i=\ell+m^*$ with $\ell,m\in L$.
If $a_j\notin L$, choose $r\in R$ with $1-ra_j\in L$ and set
$b_i=a_i$ for $i\ne j$, $b_j=a_j-r^*m$. Then
\[
 \sum_i b_i^*b_i
 =\ell+m^*(1-ra_j)-a_j^*r^*m+m^*rr^*m\in L.
\]
Weak reality gives $b_j\in L$, hence $a_j=b_j+r^*m\in L$, a
contradiction. Thus $L$ is real. By
Proposition~\ref{prop:cyclicreal}, every weakly real simple module
is real. A weakly real semisimple module has weakly real simple
summands, which are therefore real; their direct sum embeds in
their real product by Lemma~\ref{lem:apprealbasic}.
Finally, if $I=\bigcap_{j=1}^nL_j$ with each $L_j$ maximal, then
$R/I$ embeds in $\bigoplus_{j=1}^nR/L_j$ and is therefore
semisimple. Apply~(2) and Proposition~\ref{prop:cyclicreal}.
\end{proof}

If $M$ is weakly real and $R/\Ann_R(M)$ is Artinian, the latter ring
is real semisimple Artinian by Lemma~\ref{lem:apprealbasic}, so
Proposition~\ref{prop:reality-semisimplicity}(1) also gives reality of $M$.

The next proposition characterizes the existence of real maximal left
ideals over an arbitrary unital $*$-ring; for affine PI algebras,
Corollary~\ref{thm:real-max-existence} relates it to finite-dimensional
$*$-representations.

\begin{proposition}[Real simple modules and representations]
\label{prop:real-simple}
Let $R$ be a unital $*$-ring. The following are equivalent:
\begin{enumerate}
\item $R$ has a real maximal left ideal;
\item $R$ has a real simple module;
\item $R$ has a real primitive ideal;
\item $R$ admits an algebraically irreducible $*$-representation
with a positive-definite abelian-group-valued inner product.
\end{enumerate}
In~(4), the inner product is biadditive and Hermitian, satisfies
$\langle au,v\rangle=\langle u,a^*v\rangle$, and is positive definite
in the sense that $\sum_i\langle u_i,u_i\rangle=0$ forces every
$u_i=0$.
More precisely, the real primitive ideals are exactly the cores
\[
 \operatorname{core}(L)=\{a\in R:aR\subseteq L\}=\Ann_R(R/L)
\]
of real maximal left ideals $L$. Each is real, self-adjoint and
primitive, and is an intersection of real maximal left ideals.
These conditions imply that $R$ is semireal. By
Proposition~\ref{prop:reality-semisimplicity}, ``real'' may be replaced by
``weakly real'' in~(1) and~(2).
\end{proposition}
\begin{proof}
For a real simple module $S$ and $0\ne s\in S$, the left ideal
$L_s=\Ann_R(s)$ is maximal, and real by
Proposition~\ref{prop:cyclicreal}, since $R/L_s\cong S$.
Moreover, $Rs=S$ gives
\[
 \Ann_R(S)=\operatorname{core}(L_s)
 =\bigcap_{0\ne t\in S}\Ann_R(t),
\]
an intersection of real maximal left ideals.
Conversely, a real maximal left ideal $L$ gives the real simple
module $R/L$. Its annihilator is real and self-adjoint by
Lemma~\ref{lem:apprealbasic}, and primitive by definition.
A proper real left ideal excludes $-1\in\Sigma R^*R$, as in the proof
of Theorem~\ref{thm:semireal-goldie}, so these conditions imply
semireality.
For (2)$\Leftrightarrow$(4), the universal pairing on a left module
$S$ is
\[
 \langle u,v\rangle=\bar u\otimes v
 \quad\text{in }\overline S\otimes_R S,
 \qquad \langle au,v\rangle=\langle u,a^*v\rangle.
\]
The target has involution $(\bar u\otimes v)^*=\bar v\otimes u$.
Reality is exactly positive definiteness of this pairing.
Conversely, any biadditive Hermitian pairing satisfying the adjoint
identity factors through the tensor product, and positive
 definiteness makes its module real. Algebraic irreducibility is
simplicity of the underlying module.
\end{proof}

\section{PI rings}\label{sec:localization}

This section proves Theorem~\ref{thm:loc-main}. Central localization
gives prime and real prime separation for arbitrary modules, and an
embedding criterion transfers the (real) Jacobson property from rings
to finitely generated modules. The real Jacobson property can fail for
formally real affine PI algebras. We introduce positive real prime
ideals to study the existence of real maximal left ideals and give
sufficient criteria for the real Jacobson property.

\subsection{Prime and real prime separation}\label{subsec:PI-prime}

Here a \emph{PI ring} satisfies a multilinear polynomial identity
over $\mathbb Z$ with a coefficient $1$; identities arising solely
from the characteristic are excluded.

\begin{definition}[Real prime-separated modules]\label{def:real-prime-separation}
For a ring with involution, a left $R$-module is \emph{real
prime-separated} if every real submodule is an intersection of
real prime submodules. The ring $R$ is \emph{real prime-separated}
if its regular left module is.
\end{definition}

A semiprime PI ring need not be Goldie (e.g.\ $\prod_{m\ge1}\mathbb R$),
so instead of a classical quotient ring we localize at a minimal prime
of the centre.

\begin{lemma}[Central localization]\label{lem:PIloc}
Let $B$ be a nonzero semiprime PI ring with centre $C$.
There is $0\ne c\in C$ such that, for every minimal prime
$\mathfrak p$ of $C$ avoiding $c$, the set $S=C\setminus\mathfrak p$
gives a localization
\[
 A=S^{-1}B
\]
that is finite-dimensional central simple over $C_{\mathfrak p}$.
If $B$ is real, this localization is real.
\end{lemma}
\begin{proof}
By \cite[Theorem~6.1.28, Corollary~6.1.22]{Rowen}, $B$ has a PI-degree
$n$ and a central polynomial $g_n$ that is not an identity of $B$; let
$c\ne0$ be one of its values. Then $B[c^{-1}]$ is Azumaya of rank $n^2$
over $C[c^{-1}]$ \cite[Corollary~6.1.36]{Rowen}, so
$K\otimes_{C[c^{-1}]}B[c^{-1}]$ is central simple of dimension $n^2$
for every homomorphism $C[c^{-1}]\to K$ to a field
\cite[Part~B, Section~8]{DF}. The centre $C$ is reduced, since a
nilpotent $z\in C$ generates a nilpotent ideal. Since the minimal primes of the
reduced ring $C$ intersect in zero and $c\ne0$, some minimal prime
avoids $c$. Fix any such $\mathfrak p$. The reduced local ring
$C_{\mathfrak p}$ has only one prime, so it is a field. As $c\in S$,
\[
 A=S^{-1}B\cong C_{\mathfrak p}\otimes_{C[c^{-1}]}B[c^{-1}]
\]
is central simple of dimension $n^2$; the map $B\to A$ may have a kernel.

Suppose $B$ is real. If $z\in\mathfrak p$, then $uz=0$ for some
$u\notin\mathfrak p$, because $\mathfrak pC_{\mathfrak p}=0$; then
$(uz^*)^*(uz^*)=(uz)(u^*z^*)=0$, so $uz^*=0$ and $z^*\in\mathfrak p$.
Thus $S^*=S$ and the involution extends to $A$. If
$\sum_i(b_i/s)^*(b_i/s)=0$, then $t\sum_ib_i^*b_i=0$ for some $t\in S$,
so $\sum_i(tb_i)^*(tb_i)=0$, every $tb_i=0$, and every $b_i/s=0$.
\end{proof}

\begin{theorem}[Prime separation for PI rings]\label{thm:PI}\label{cor:PIbounded}
Let $R$ be a PI ring. Every semiprime submodule of every $R$-module
is an intersection of prime submodules. If $R$ has an involution,
an $R$-module is real if and only if it is semiprime and weakly real,
and every real submodule is an intersection of real prime submodules.
In particular, weakly real left ideals are real. No finite generation
is assumed.
\end{theorem}
\begin{proof}
Let $M$ be semiprime and $0\ne q\in M$. Put
$L=Rq$, $T=\Ann_R(L)$ and $B=R/T$. Then $B$ is nonzero and semiprime,
$L$ is faithful over $B$, and $L\cap TM=0$ by
Lemma~\ref{lem:prime-semiprime-annihilators}(2) and semiprimeness of
$M$. Choose $A=S^{-1}B$ as in Lemma~\ref{lem:PIloc} and put
$V=S^{-1}(M/TM)$. The map $M\to V$ does not kill $q$: otherwise
$\tilde sq\in L\cap TM=0$ for a lift $\tilde s$ of some $s\in S$; since
$s$ is central in $B$, this gives $\tilde s(rq)=r\tilde sq=0$ for all
$r$, so $sL=0$, contrary to faithfulness.

For every nonzero $B$-submodule $E\subseteq V$, the $A$-module
$S^{-1}E$ is nonzero, so its annihilator in the simple ring $A$ is
zero. Hence, by centrality, $bE=0$ forces $b/1=0$; conversely
$\ker(B\to A)$ kills $V$. Thus $\Ann_B(E)=\ker(B\to A)$ for all such $E$, so $V$ is prime over
$B$, hence over $R$, and $P=\ker(M\to V)$ is a prime submodule missing
$q$.

If $M$ is also weakly real, then $T$ is real and self-adjoint by
Lemma~\ref{lem:apprealbasic}, so $A$ is real, and $V$ is real over $R$
by Proposition~\ref{prop:reality-semisimplicity}(1) and restriction of
scalars (Lemma~\ref{lem:apprealbasic}). Hence $M/P$, which
embeds in $V$, is real. These real prime quotients embed $M$ in their
product, so $M$ is real; the converse is
Lemma~\ref{lem:apprealbasic}. Applying the construction to real
quotients proves real prime separation, and a weakly real left ideal
is real because its quotient is semiprime and weakly real
(Proposition~\ref{prop:cyclicreal}).
\end{proof}

\begin{remark}\label{rem:PIcompare}
Theorem~\ref{thm:PI} applies to rings module-finite over their centre,
for instance matrix rings over commutative rings. For commutative rings
with trivial involution, the characterization of real modules is the
annihilator criterion of \cite[Proposition~4.3]{Cdec}; semiprimeness
cannot be omitted \cite[Example~4.4]{Cdec}, \cite[Theorem~11]{JS}.
For left ideals, the theorem extends \cite[Corollary~1]{C22b} to all
PI $*$-rings and answers the commutative nontrivial-involution question
in \cite{Cdec}. The prime-separation argument uses only simplicity
of $A$, so it applies whenever every nonzero semiprime factor ring has a
nonzero simple central localization, e.g.\ to $E\otimes_kC$ with $E$
simple with centre $k$ and $C$ commutative
\cite[Introduction and Corollary~3.12]{Gogic}.
\end{remark}

\subsection{Jacobson rings and modules}\label{subsec:PI-maximal}

\begin{definition}[Jacobson modules]\label{def:maximal-separation}
A left $R$-module $M$ is \emph{Jacobson} if every prime submodule
is an intersection of maximal submodules. If $R$ has an involution,
$M$ is \emph{real Jacobson} if every real prime submodule is an
intersection of real maximal submodules.
Thus $R$ is left Jacobson, respectively real left Jacobson, exactly
when its regular left module has the corresponding property.
\end{definition}

For commutative rings, Jacobson modules in this sense are the
\emph{Hilbert modules} of
\cite{MSS}; our left Jacobson rings are the weakly left Jacobson rings
of \cite{CS}, and their strongly left Jacobson rings are the rings that
are both prime-separated and left Jacobson. For a prime-separated module, the Jacobson property means
that every semiprime submodule is an intersection of maximal ones,
and similarly in the real case.

\begin{definition}[Prime embedding property]\label{def:prime-embedding}
A ring $R$ has the \emph{prime embedding property} if every finitely
generated prime left $R$-module $V$ embeds $R$-linearly in a product
of copies of $R/\Ann_R(V)$. A ring with involution has the \emph{real
prime embedding property} if this holds for every finitely generated
real prime left module.
\end{definition}

\begin{theorem}[Jacobson transfer]\label{thm:jacobson-transfer}
Let $R$ be a ring.
\begin{enumerate}
\item If $R$ has the prime embedding property, the following are
equivalent:
\begin{enumerate}
\item[(a)] $R$ is Jacobson.
\item[(b)] $R$ is left Jacobson.
\item[(c)] Every finitely generated left $R$-module is Jacobson.
\end{enumerate}
\item If $R$ has an involution and the real prime embedding property,
the following are equivalent:
\begin{enumerate}
\item[(a)] $R$ is real Jacobson.
\item[(b)] $R$ is real left Jacobson.
\item[(c)] Every finitely generated left $R$-module is real Jacobson.
\end{enumerate}
\end{enumerate}
The implications \textup{(c)}$\Rightarrow$\textup{(b)}$\Rightarrow$\textup{(a)}
hold without the embedding hypotheses.
\end{theorem}
\begin{proof}
We prove both parts together, imposing reality in~(2).
(c)$\Rightarrow$(b) is the case $M=R$. For (b)$\Rightarrow$(a), regard a prime two-sided ideal $T$, real
in~(2), as a prime left ideal. If $T=\bigcap_i L_i$ with maximal
left ideals $L_i$, real in~(2), then
\[
 T\subseteq\Ann_R(R/L_i)\subseteq L_i,\qquad
 T=\bigcap_i\Ann_R(R/L_i).
\]
These annihilators are primitive, and real primitive in~(2).

For (a)$\Rightarrow$(c), let $P$ be a prime submodule of a finitely
generated module $M$, real in~(2), and let $a\notin P$. Put
$V=M/P$ and $T=\Ann_R(V)$. Then $T$ is prime, and real in~(2)
by Lemma~\ref{lem:apprealbasic}. By (a), $T$ is an intersection of
annihilators of simple modules, real simple modules in~(2).
For each such $U$, $\Ann_R(U)=\bigcap_{0\ne u\in U}\Ann_R(u)$ with
$R/\Ann_R(u)\cong U$, so the maximal left ideals containing $T$, real
maximal in~(2), intersect in $T$.

The embedding property supplies an $R$-linear map $f:V\to R/T$
with $f(a+P)\ne0$. Choose a maximal left ideal $L\supseteq T$,
real in~(2), missing a representative of $f(a+P)$. The composite
\[
 M\twoheadrightarrow V\xrightarrow{\,f\,}R/T\twoheadrightarrow R/L
\]
is nonzero at $a$, hence surjective because $R/L$ is simple.
Its kernel is maximal, contains $P$ and misses $a$; in~(2) it is
real maximal because $R/L$ is real. Varying $a$ proves (c).
\end{proof}

The first Weyl algebra in characteristic zero is Jacobson but not
left Jacobson \cite[Theorem~7]{CS}. The real converse can also fail:
in separate work we will exhibit a formally real division $*$-ring
$D$ for which $D[z]$ is real Jacobson but not real left Jacobson.

We now verify the prime embedding property for all PI rings.

\begin{lemma}[Prime embedding for PI rings]\label{lem:PIembedding}
Every PI ring $R$ has the prime embedding property and, if it has an
involution, the real prime embedding property. In fact, every finitely
generated prime left $R$-module $V$ embeds $R$-linearly in
$(R/\Ann_R(V))^r$ for some $r$.
\end{lemma}
\begin{proof}
Let $V$ be such a module and put $B=R/\Ann_R(V)$. Then $B$ is a
prime PI ring and $V$ is a finitely generated faithful prime
$B$-module. By Posner's theorem \cite{Posner},
\cite[Part~B, Section~6]{DF}, $Q=S^{-1}B$ is simple Artinian for
$S=Z(B)\setminus\{0\}$. An equality $sv=0$ with $s\in S$ and
$v\ne0$ would give $sBv=0$, contradicting $\Ann_B(Bv)=0$.
Hence $V\hookrightarrow S^{-1}V$.
A finite generating surjection $Q^r\twoheadrightarrow S^{-1}V$
splits, giving an embedding $f:V\to Q^r$. Choose a common central
denominator $s$ for the images of finitely many generators of $V$.
Then $v\mapsto sf(v)$ embeds $V$ in $B^r$.
The real prime embedding property follows because every real prime
module is prime.
\end{proof}

Theorem~\ref{thm:jacobson-transfer} now gives the Jacobson
equivalences in Theorem~\ref{thm:loc-main}.

\begin{corollary}[Amitsur--Procesi for modules]\label{cor:affine-PI-maximal}
Let $R$ be an affine PI algebra over a field $k$. In every finitely
generated left $R$-module, every semiprime submodule is an intersection
of maximal submodules, and every maximal submodule has finite
$k$-codimension. If $R$ has an involution and is real Jacobson, the same
assertion holds with real submodules and real maximal submodules.
\end{corollary}
\begin{proof}
The algebra $R$ is Jacobson by Amitsur--Procesi
\cite[Part~B, Theorem~11.7]{DF}, so
Theorems~\ref{thm:PI} and~\ref{thm:jacobson-transfer}, together with
Lemma~\ref{lem:PIembedding}, give the intersection assertions.
Let $L$ be a maximal submodule of a finitely generated left $R$-module
$M$, and put $U=M/L$. Kaplansky's theorem makes
$A=R/\Ann_R(U)$ finite-dimensional over its centre $F$.
Since $A$ is affine, Artin--Tate and Zariski's lemma give
$[F:k]<\infty$. Thus $A$ and its cyclic module $U$ are
finite-dimensional over $k$, so $L$ has finite $k$-codimension;
see also \cite[Part~B, Theorem~11.5]{DF}.
\end{proof}

\begin{remark}[A Jacobson centre is not enough]\label{rem:jacobson-centre}
For a field $k$ and $A=k[[t]]$, the ring
$R=\left(\begin{smallmatrix}k&A\\0&A\end{smallmatrix}\right)$
satisfies $[x_1,x_2][x_3,x_4]=0$ (commutators are strictly upper
triangular) and has centre $kI_2$, but it maps onto $A$, whose zero
ideal is prime while $(t)$ is its only primitive ideal; so $R$ is not
Jacobson. It is right Noetherian,
so it satisfies hypothesis~(i) of Theorem~\ref{thm:main}, but not~(ii):
$Re_{11}R\cong k\oplus A$, with $R$ acting on the left through
$R\to k$, is not finitely generated as a left ideal.
\end{remark}

\begin{proposition}[Real primitive ideals of PI rings]
\label{prop:PI-real-primitive}
Let $R$ be a PI ring with involution. A two-sided ideal $J$ of $R$ is
real primitive if and only if it is real and primitive. If $R$ is an
affine PI algebra over a field $k$, every primitive ideal of $R$ has
finite $k$-codimension.
\end{proposition}
\begin{proof}
Real primitive ideals are real and primitive by
Proposition~\ref{prop:real-simple}. Conversely, if $J$ is real and
primitive, then $J=J^*$ (Subsection~\ref{subsec:realideals}), and the
primitive PI factor $R/J$ is formally real, and simple Artinian by
Kaplansky's theorem \cite[Part~B, Theorem~4.13]{DF}. Its
simple modules are real by Proposition~\ref{prop:reality-semisimplicity}(1) and
faithful, so $J$ is real primitive. The last assertion holds because
simple modules over affine PI algebras are finite-dimensional over $k$
(proof of Corollary~\ref{cor:affine-PI-maximal}).
\end{proof}

\subsection{A formally real affine PI counterexample}
\label{sec:PIcounter}

An affine PI algebra need not be real Jacobson, even though it is
real prime-separated by Theorem~\ref{thm:PI}.

\begin{example}[A twisted adjoint]\label{ex:twistedreal}
Put
\[
 C=\mathbb R[x^{\pm1},y^{\pm1}],\qquad B=M_4(C),\qquad
 S=\operatorname{diag}(1,x,y,-xy),
\]
and define $a^*=S^{-1}a^{\mathsf T}S$. Then $B$ is an affine Noetherian
PI $\mathbb R$-algebra, finite over its centre. Its involution is
formally real, but $B$ has no real primitive ideals and no real
maximal left ideals. Thus the real prime ideal $0$ is not an
intersection of real primitive ideals, and $B$ is not real Jacobson.
\end{example}
\begin{proof}
Since $S$ is symmetric and invertible over $C$, the formula defines an
$\mathbb R$-linear involution. The algebra is finite free over the affine
commutative domain $C$, and is prime, Noetherian and PI.

\emph{Formal reality.} We first observe that
\[
 p_1+xp_2+yp_3-xyp_4=0,
 \qquad p_k\text{ sums of squares in }C,
\]
forces every $p_k=0$. Order Laurent-monomial exponents in $\mathbb Z^2$
lexicographically. A nonzero sum of squares has a smallest exponent in
$2\mathbb Z^2$ with positive coefficient: the smallest monomial of a
square has twice the exponent of the smallest monomial of its square
root, and these positive coefficients cannot cancel in a sum.
The nonzero summands in the displayed equality therefore have smallest
exponents in distinct parity classes $(0,0),(1,0),(0,1),(1,1)$.
One has a uniquely smallest exponent, whose coefficient cannot cancel.
This proves the observation.

Write $s_1=1,s_2=x,s_3=y,s_4=-xy$. If $\sum_\nu a_\nu^*a_\nu=0$,
the $(j,j)$ entry, multiplied by $s_j$, is
\[
 \sum_{k=1}^4 s_k\left(\sum_\nu(a_\nu)_{kj}^{\,2}\right)=0.
\]
The observation forces all four sums of squares to vanish, hence every
entry in column $j$ of every $a_\nu$ is zero. This holds for every $j$,
so $B$ is formally real. As $B$ is prime, so is its real left ideal $0$.

\emph{Absence of real primitive ideals and real maximal left ideals.} Every primitive
ideal of $B$ is $M_4(\mathfrak m)$ for a maximal ideal $\mathfrak m$
of $C$. A real maximal left ideal would give a real primitive ideal
by Proposition~\ref{prop:real-simple}, so it suffices to show that none of these primitive
quotients is formally real. The residue field is $\mathbb R$ or $\mathbb C$.
In the complex case the induced involution fixes the centre pointwise;
a central $i$ then gives $1^*1+i^*i=0$, contradicting reality.
In the real case, $\mathfrak m$ is evaluation at
$(a,b)\in(\mathbb R^\times)^2$, and the quotient involution is adjoint to
\[
 \langle1,a,b,-ab\rangle.
\]
This form is always indefinite: positive definiteness would require
$a,b,-ab>0$, and negative definiteness is excluded by the coefficient
$1$. Choose a nonzero isotropic vector $z$. A matrix $h$ with first
column $z$ and all other columns zero satisfies $h\ne0$ and $h^*h=0$.
Thus this quotient is not formally real either.
\end{proof}

\subsection{Positive involutions and \texorpdfstring{$*$}{*}-representations}
\label{subsec:star-representations}

Let $A$ be a finite-dimensional simple algebra with
involution $*$, centre $K$ and fixed centre $F=\{a\in K:a^*=a\}$. Since
$\Trd_{A/K}(a^*)=\Trd_{A/K}(a)^*$, the reduced trace of a hermitian
square lies in $F$, and
\[
 T_*(x)=\Trd_{A/K}(x^*x)\qquad(x\in A)
\]
is a quadratic form over $F$, the \emph{involution trace form}. The
involution is \emph{positive at an ordering} of $F$ if $T_*$ is positive
definite \cite[Definition~1.1]{PS}. Positivity implies formal reality:
applying $\Trd_{A/K}$ to $\sum_ix_i^*x_i=0$ forces every $x_i=0$.

\begin{remark}[Normalizations of the trace form]\label{rem:trace-normalization}
Our normalization agrees with Procesi--Schacher \cite[Definition~1.1]{PS}
and Klep--Unger \cite[Section~2]{KU} for involutions of the first kind.
In the unitary case they use $2T_*$, which gives the same notion of
positivity.
\end{remark}

For $k=\mathbb R$, the embedding in the next lemma follows from
\cite[Proposition~2.7]{CHMN15}.

\begin{lemma}[Finite-dimensional real algebras over a real closed field]
\label{lem:realclosed-positive}
Let $k$ be real closed and $A$ a nonzero finite-dimensional $k$-algebra
with a $k$-linear involution.
\begin{enumerate}
\item If $A$ is real, it admits a unital $*$-embedding into $M_N(k)$
with transpose involution.
\item If $A$ is simple, then $A$ is real if and only if its fixed
centre is $k$ and its involution is positive at the unique ordering of $k$.
\end{enumerate}
\end{lemma}
\begin{proof}
(1) Suppose $A$ is real. Then $A$ is semisimple, and the involution
fixes the central
idempotent $e$ of each simple factor (else $e^*e=0$). Embedding the
factors block-diagonally, we may take $A$ simple. Its centre $K$ is
$k$ or $k(i)$; in the latter case reality forces conjugation on $K$,
since a symmetric $i$ would give $1^*1+i^*i=0$. Thus the fixed centre
is $k$.

For $0\ne x\in A$, set $s=x^*x$. The algebra $k[s]$ is real with
trivial involution, hence reduced. Being finite-dimensional, it is a
product of fields, which by real closedness are copies of $k$ or
$k(i)$; a $k(i)$ factor would give symmetric $e\ne0,u$ with $e^2=e$,
$u^2=-e$ and $e^*e+u^*u=0$, contrary to reality. Thus
$s=\sum_j\lambda_je_j$ with $\lambda_j\in k$ and nonzero orthogonal
symmetric coordinate idempotents $e_j$. If $\lambda_j<0$, then
$\mu=\sqrt{-\lambda_j}\in k$ gives
$(xe_j)^*(xe_j)+(\mu e_j)^*(\mu e_j)=0$, again contradicting reality.
Hence every $\lambda_j\ge0$, and at least one is positive since
reality gives $s\ne0$. Each $\Trd_{A/K}(e_j)$ is a positive integer,
the rank of the idempotent after extension to a splitting field. Therefore
\[
 T_*(x)=\sum_j\lambda_j\Trd_{A/K}(e_j)>0.
\]
The polarization
$\langle x,y\rangle=\frac12\Trd_{A/K}(x^*y+y^*x)$ is a positive definite
symmetric $k$-bilinear form with
$\langle ax,y\rangle=\langle x,a^*y\rangle$. In an orthonormal basis,
which exists because positive elements of $k$ are squares, the faithful
left regular action gives a unital $*$-embedding
$A\hookrightarrow M_{\dim_k A}(k)$ with transpose.

(2) For simple real $A$, the argument above identifies its fixed centre
with $k$ and proves positivity. Conversely, applying the reduced trace
to $\sum_ix_i^*x_i=0$ and using positivity forces every $x_i=0$.
\end{proof}

Let $R$ be a $k$-algebra with a $k$-linear involution and $E$ an ordered
field extension of $k$. A \emph{$*$-representation} of $R$ over $E$ is a
unital $k$-algebra homomorphism $\pi\colon R\to\operatorname{End}_E(V)$,
where $V$ is a nonzero finite-dimensional $E$-space with a positive
definite symmetric bilinear form such that
\[
 \langle\pi(a)u,v\rangle=\langle u,\pi(a^*)v\rangle\qquad(a\in R).
\]
It is \emph{irreducible} if $V$ has no proper nonzero invariant
$E$-subspace. Extending scalars to a real closure of $E$ preserves positive
definiteness, and over a real closed field an orthonormal basis turns
$\pi$ into a unital $*$-homomorphism $R\to(M_n(E),{}^{\mathsf T})$.
Conversely, such a homomorphism over a formally real field $E$ is a
$*$-representation on $E^n$ with the standard form, for any ordering of
$E$. The following decomposition is also used for free algebras in
Section~\ref{sec:freereal}.

\begin{lemma}[Orthogonal decomposition]\label{lem:star-decomposition}
Let $\pi$ be a $*$-representation of $R$ over an ordered field extension
$E$ of $k$.
\begin{enumerate}
\item $\pi$ is an orthogonal direct sum of irreducible
$*$-representations, and $\ker\pi$ is real.
\item If $\pi$ is irreducible, then $\ker\pi$ is prime and $E\pi(R)$ is
simple Artinian. If moreover $R$ is affine and $E$ is algebraic over
$k$, then $\ker\pi$ is real primitive of finite $k$-codimension.
\end{enumerate}
The decomposition in~\textup{(1)} also holds over $\mathbb C$ for
positive definite Hermitian forms, with the same proof.
\end{lemma}
\begin{proof}
(1) If $W\subseteq V$ is invariant, the adjoint identity makes $W^\perp$
invariant, and positive definiteness gives $V=W\oplus W^\perp$.
Induction on $\dim_EV$ gives the decomposition. If
$\sum_ia_i^*a_i\in\ker\pi$, then
$\sum_i\langle\pi(a_i)v,\pi(a_i)v\rangle=0$ for every $v\in V$, so every
$\pi(a_i)=0$; thus $\ker\pi$ is real.

(2) If two-sided ideals $I,J$ satisfy $IJ\subseteq\ker\pi$ and
$J\nsubseteq\ker\pi$, then $\pi(J)V$ is a nonzero invariant subspace,
hence equals $V$, and $\pi(I)V=\pi(IJ)V=0$. Thus $\ker\pi$ is prime.
An $E\pi(R)$-submodule of $V$ is precisely a $\pi(R)$-invariant
$E$-subspace, so irreducibility makes $V$ simple over $E\pi(R)$.
This action is faithful because $E\pi(R)\subseteq\operatorname{End}_E(V)$;
hence the finite-dimensional algebra $E\pi(R)$ is simple Artinian.
If $R$ is affine and $E$ is algebraic over
$k$, the matrices of finitely many generators, in any basis of $V$,
have entries in a finite extension of $k$. Hence $\pi(R)\cong R/\ker\pi$
is finite-dimensional over $k$ and prime, hence simple Artinian, and it
is formally real by~(1). Its simple modules are real by
Proposition~\ref{prop:reality-semisimplicity}(1) and faithful, so
$\ker\pi$ is real primitive.
\end{proof}

Lemma~\ref{lem:transpose-specialization} will be used in
Corollaries~\ref{cor:positive-real-primitive}--\ref{thm:real-max-existence}
and Theorem~\ref{thm:realPI}.

\begin{lemma}[Specialization]\label{lem:transpose-specialization}
Let $k$ be a formally real field and $R$ an affine $k$-algebra with a
$k$-linear involution. Let $\varphi$ be a $*$-representation of $R$
over an ordered field extension $E$ of $k$, and let $b\in R$ with
$\varphi(b)\ne0$. Then there are a real closure $K$ of $k$ and a
$*$-representation $\psi$ of $R$ over $K$ with
$\ker\varphi\subseteq\ker\psi$ and $\psi(b)\ne0$.
\end{lemma}
\begin{proof}
Extending scalars to a real closure $\Omega$ of $E$ and choosing an
orthonormal basis, we may assume that
$\varphi(R)\subseteq(M_n(\Omega),{}^{\mathsf T})$. Let $K$ be the
subfield of $\Omega$ consisting of elements algebraic over $k$; it is
a real closure of $k$. Choose finitely many $k$-algebra generators of
$R$, and let $C\subseteq\Omega$ be the $k$-algebra generated
by all entries of their images. Then $\varphi(R)\subseteq M_n(C)$, and
$\varphi(b)$ has a nonzero entry $c\in C$. Write
$C=k[t_1,\ldots,t_m]/(f_1,\ldots,f_r)$ with $c$ represented by $g(t)$.
The system $f_1=\cdots=f_r=0$, $g(t)z=1$ has a solution in $\Omega$,
hence, by real closed field transfer \cite[Chapter~5]{BCR}, one in
$K$. Evaluation at it gives a $k$-algebra homomorphism
$\varepsilon\colon C\to K$ with $\varepsilon(c)\ne0$. Applied
entrywise, $\varepsilon$ commutes with transposition, so
$\psi=\varepsilon\circ\varphi$ is a $*$-representation over $K$
with $\ker\varphi\subseteq\ker\psi$ and $\psi(b)\ne0$.
\end{proof}

\subsection{Positive real primes and the real Jacobson property}
\label{subsec:PI-real-primitive}\label{subsec:positivity-criterion}

We now ask when real maximal left ideals exist and when $R$ is real
Jacobson. Throughout this subsection, $k$ is a formally real field and
$R$ is an affine PI $k$-algebra with a $k$-linear involution.
For a prime PI ring $B$, the \emph{central simple ring of quotients} is
$(Z(B)\setminus0)^{-1}B$ \cite{Posner}; an involution of $B$ extends to
it. A real prime two-sided ideal $T$ of $R$ is \emph{positive} if the
involution induced on the central simple ring of quotients of $R/T$
is positive at some ordering of its fixed centre.

We first describe positive real primes by $*$-representations
(Proposition~\ref{prop:positive-prime-kernels}) and specialize this to
real primitive ideals (Corollary~\ref{cor:positive-real-primitive}).
The main result relates positivity of all real primes to the real
Jacobson property (Theorem~\ref{thm:realPI}).

\begin{proposition}[Positive real primes and representations]
\label{prop:positive-prime-kernels}
For a real prime two-sided ideal $T$ of $R$, the following conditions
are equivalent:
\begin{enumerate}
\item $T$ is positive;
\item $T$ is the kernel of a $*$-representation over some ordered field
extension of $k$;
\item $T$ is the kernel of an irreducible $*$-representation over some
ordered field extension of $k$.
\end{enumerate}
In~\textup{(2)} and~\textup{(3)}, the ordered field can be taken to be
the fixed centre of the central simple ring of quotients of $R/T$, with
an ordering at which its involution is positive.
\end{proposition}
\begin{proof}
Put $B=R/T$, and let $Q$ be its central simple ring of quotients and $F$ the
fixed centre of $Q$. Since the involution is $k$-linear, $k\subseteq F$.

(1)$\Rightarrow$(3). Fix an ordering of $F$ at which the involution
trace form $T_*$ of $Q$ is positive definite. Its polarization
\[
 \langle x,y\rangle=\tfrac12\Trd_{Q/Z(Q)}(x^*y+y^*x)
\]
is a positive definite symmetric $F$-bilinear form on $Q$, and
$\langle ax,y\rangle=\langle x,a^*y\rangle$ for $a,x,y\in Q$. Moreover
$FB=Q$: a nonzero $s\in Z(B)$ satisfies $s^{-1}=s^*(s^*s)^{-1}$ with
$s^*s\in F$. Let $L$ be a simple left ideal of $Q$. An $F$-subspace of
$L$ invariant under $B$ is invariant under $FB=Q$, so left
multiplication $R\to B\to\operatorname{End}_F(L)$, with the restricted
form, is an irreducible $*$-representation over the ordered field $F$.
Its kernel is $T$, because the simple $Q$-module $L$ is faithful.

(3)$\Rightarrow$(2) is trivial.

(2)$\Rightarrow$(3) and (2)$\Rightarrow$(1). Let $\varphi$ be a
$*$-representation with kernel $T$. Since extending scalars does not
change the kernel, we may assume that it is defined over a real closed
field $\Omega$. By
Lemma~\ref{lem:star-decomposition}, $T$ is the intersection of the
finitely many kernels of the irreducible summands of $\varphi$; since
$T$ is prime, one summand $\psi$ has kernel $T$. This proves~(3).
Identifying $B$ with $\psi(R)$, the algebra $A=\Omega B$ is simple
Artinian by Lemma~\ref{lem:star-decomposition}(2), and formally real as
a $*$-subalgebra of $(M_n(\Omega),{}^{\mathsf T})$ in an orthonormal
basis. Elements of $Z(B)$ commute with $\Omega B=A$, so $Z(B)$ lies in
the field $Z(A)$. Every nonzero $s\in Z(B)$ is therefore invertible in
$A$. Thus the inclusion $B\hookrightarrow A$ extends to an embedding
$Q=(Z(B)\setminus\{0\})^{-1}B\hookrightarrow A$ by
$b/s\mapsto bs^{-1}$; injectivity follows since $bs^{-1}=0$ implies
$b=0$. Since $A$ is the $\Omega$-span of $B$ and $\Omega\subseteq Z(A)$,
multiplication gives a surjection
\[
 Q\otimes_{Z(Q)}Z(A)\longrightarrow A.
\]
Its domain is central simple over $Z(A)$, so it is an isomorphism, and
reduced traces of $Q$ are restrictions of those of $A$. By
Lemma~\ref{lem:realclosed-positive}(2), the fixed centre of $A$ is $\Omega$
and $\Trd_{A/Z(A)}(a^*a)>0$ for $a\ne0$. Hence $F\subseteq\Omega$
inherits an ordering at which $\Trd_{Q/Z(Q)}(q^*q)>0$ for every $q\ne0$,
so the involution of $Q$ is positive.
\end{proof}

For primitive ideals, maximality allows specialization to a real
closure without changing the kernel.

\begin{corollary}[Positive real primitive ideals]
\label{cor:positive-real-primitive}
The positive real primitive ideals of $R$ are precisely the kernels of
irreducible $*$-representations over real closures of $k$. If $k$ is
real closed, every real primitive ideal is positive.
\end{corollary}
\begin{proof}
Let $J$ be a real primitive ideal, hence a maximal two-sided ideal.
If $J$ is positive, it is the kernel of
a $*$-representation by Proposition~\ref{prop:positive-prime-kernels}, and
Lemma~\ref{lem:transpose-specialization} with $b=1$ gives a
$*$-representation $\psi$ over a real closure of $k$ with
$J\subseteq\ker\psi$. The kernels of the irreducible summands of $\psi$
are proper and contain $J$, so they equal $J$. Conversely, the kernel of
an irreducible $*$-representation over a real closure of $k$ is real
primitive by Lemma~\ref{lem:star-decomposition}(2) and positive by
Proposition~\ref{prop:positive-prime-kernels}. If $k$ is real closed, then $R/J$ is
finite-dimensional, simple and formally real by
Proposition~\ref{prop:PI-real-primitive}, and
Lemma~\ref{lem:realclosed-positive}(1) gives a faithful $*$-representation
of $R/J$ over $k$; thus $J$ is positive by
Proposition~\ref{prop:positive-prime-kernels}.
\end{proof}

\begin{corollary}[Existence]\label{thm:real-max-existence}
$R$ has a $*$-representation over some ordered field extension
of $k$ if and only if it has one over a real closure of $k$, if and only
if some real prime ideal of $R$ is positive. In this case $R$ has a real
maximal left ideal; if $k$ is real closed, the converse holds.
\end{corollary}
\begin{proof}
Lemma~\ref{lem:transpose-specialization} with $b=1$ gives the first
equivalence. By Lemma~\ref{lem:star-decomposition} and
Corollary~\ref{cor:positive-real-primitive}, a representation over a
real closure has an irreducible summand with positive real primitive
kernel, hence yields a real maximal left ideal by
Proposition~\ref{prop:real-simple}. Conversely, a positive real prime
gives a $*$-representation by Proposition~\ref{prop:positive-prime-kernels}.
When $k$ is real closed, Corollary~\ref{cor:positive-real-primitive} and
Proposition~\ref{prop:real-simple} give the final converse.
\end{proof}

We now pass from individual kernels to intersections. The main step
is to combine a family of representations whose kernels have prime
intersection into a single representation with the same kernel.

\begin{theorem}[Positive real primes and the real Jacobson property]
\label{thm:positive-prime}\label{thm:realPI}
Let $k$ be a formally real field and $R$ an affine PI $k$-algebra with a
$k$-linear involution.
\begin{enumerate}
\item[(i)] A real prime two-sided ideal $T$ of $R$ is positive if and
only if it is an intersection of positive real primitive ideals.
\item[(ii)] Every real prime two-sided ideal of $R$ is positive if and
only if $R$ is real Jacobson and every real primitive ideal of $R$ is
positive.
\end{enumerate}
\end{theorem}
\begin{proof}
(i) Suppose $T$ is positive. Proposition~\ref{prop:positive-prime-kernels}
gives a $*$-representation with kernel $T$. For $b\notin T$,
Lemma~\ref{lem:transpose-specialization} gives a $*$-representation
over a real closure of $k$ whose kernel contains $T$ and misses $b$; by
Lemma~\ref{lem:star-decomposition}, one of its irreducible summands
still misses $b$. Its kernel is positive real primitive by
Corollary~\ref{cor:positive-real-primitive}. Thus $T$ is the intersection
of these kernels.

Conversely, by Corollary~\ref{cor:positive-real-primitive}, write $T$ as
an intersection of kernels of irreducible $*$-representations over
real closures of $k$. Fix a multilinear
identity of $R$ of degree $d$. For an irreducible $*$-representation
$\pi$ over a real closure $L$ of $k$, the algebra $A=L\pi(R)$ is simple
by Lemma~\ref{lem:star-decomposition}(2), and satisfies the same
identity by multilinearity. Since $Z(A)$ is $L$ or $L(i)$, the
PI-degree bound gives $\dim_LA\le2(d/2)^2\le d^2$
\cite[Part~B, Lemma~4.14(a),(c)]{DF}. The representation space is cyclic over
$A$, so its dimension is at most $d^2$.

Choose a common multiple $N$ of $1,\ldots,d^2$. Taking orthogonal sums
of copies makes all representations in such a family have
dimension $N$ without changing their kernels. In orthonormal bases,
write them as $\varphi\colon R\to(M_N(L_\varphi),{}^{\mathsf T})$.
Let $C\subseteq\prod_\varphi L_\varphi$ be generated over $k$ by the
entry tuples of the images of finitely many generators of $R$. Then
$R/T\hookrightarrow M_N(C)$ compatibly with transposition, and $C$ is
Noetherian, reduced and real.

Let $S$ be the set of nonzerodivisors of $C$. The total ring of fractions
$S^{-1}C$ is a finite product of fields $E_1,\ldots,E_s$
\cite[Tags~00EU, 00EW and 02LX]{Stacks}, and $C\hookrightarrow S^{-1}C$. This ring
remains real: a relation $\sum_j(b_j/u)^2=0$, with $u\in S$, gives
$\sum_jb_j^2=0$ in $C$, hence every $b_j=0$. Reality passes to the
direct factors, so each $E_i$ is formally real.

Choose an ordering of each $E_i$. The composites
$R\to M_N(C)\to M_N(E_i)$ are $*$-representations with the standard
forms, and their kernels intersect in $T$ because
$C\hookrightarrow\prod_iE_i$. The product of these finitely many
kernels is contained in $T$, so primeness forces one of them to equal $T$.
Proposition~\ref{prop:positive-prime-kernels} now makes $T$ positive.

(ii) Real primitive ideals are real prime. If every real prime ideal is
positive, part~\textup{(i)} writes each of them as an
intersection of real primitive ideals, so $R$ is real Jacobson.
Conversely, if $R$ is real Jacobson, every real prime ideal is an
intersection of real primitive ideals; if these are positive, so is the
real prime ideal, by part~\textup{(i)}.
\end{proof}

Positivity is hard to verify directly; for Azumaya algebras, a condition
on reduced traces of hermitian squares suffices.

\begin{corollary}[Procesi--Schacher condition]
\label{cor:azumaya-trace-square}
Let $R$ be an Azumaya algebra whose centre $C$ is an affine algebra
over a formally real field $k$, and let $*$ be a $k$-linear involution
on $R$.
If $\Trd_{R/C}(a^*a)\,1_R\in\Sigma R^*R$ for every $a\in R$, then every
real prime two-sided ideal of $R$ is positive; in particular, $R$ is real
Jacobson by Theorem~\ref{thm:realPI}(ii).
\end{corollary}
\begin{proof}
Since $R$ is finite over $C$, it is an affine PI $k$-algebra.
Let $T$ be a real prime two-sided ideal and put $\mathfrak p=T\cap C$.
By the Azumaya ideal correspondence, $T=\mathfrak pR$, so the central simple ring of quotients of $R/T$ is $Q=K\otimes_{C/\mathfrak p}R/T$ with
$K=\operatorname{Frac}(C/\mathfrak p)$; let $F$ be its fixed centre.
Reduced trace commutes with this base change. Elements of $Q$ have the
form $a/s$ with $0\ne s\in C/\mathfrak p$; clearing such denominators
shows that $Q$ is real, and
$\operatorname{Trd}((a/s)^*(a/s))=\operatorname{Trd}(a^*a)/(s^*s)$
carries the trace-square identities to $Q$, dividing their witnesses
by $s$. Thus $P=F\cap\Sigma Q^*Q$
contains the values $\operatorname{Trd}(q^*q)$ of the involution trace
form. It is a proper
preordering of $F$: it contains squares, is closed under addition, and
reality excludes $-1$. For multiplication, if
$\alpha=\sum_ia_i^*a_i$ and $\beta=\sum_jb_j^*b_j$ belong to $P$,
centrality gives
$\alpha\beta=\sum_jb_j^*\alpha b_j=\sum_{i,j}(a_ib_j)^*(a_ib_j)\in P$. An ordering containing $P$
makes this form nonnegative, hence positive definite because it is
nonsingular. Thus $T$ is positive.
\end{proof}

\begin{remark}[Trace squares, positivity and real Jacobsonness]
\label{rem:positivity-not-necessary}
By Corollary~\ref{cor:azumaya-trace-square} and
Theorem~\ref{thm:realPI}(ii), the trace-square condition implies
positivity of every real prime ideal, which in turn implies real
Jacobsonness. Neither implication can be
reversed, even for central simple algebras.

For central simple algebras, the trace-square condition is equivalent to
the Procesi--Schacher conjecture \cite[Remark~2.5]{KU}. Klep and Unger
give a split degree-three counterexample whose involution is positive
at some ordering but whose trace-square condition fails
\cite[Theorem~3.2]{KU}. Since zero is its only proper prime ideal, this
shows that the first implication cannot be reversed.

For the second implication, let $k=\mathbb R(x,y)$, and let $A=M_4(k)$
carry the involution adjoint to the form $\langle1,x,y,-xy\rangle$.
This algebra is simple and formally real (clear denominators in
Example~\ref{ex:twistedreal}).
It is real Jacobson and has real maximal left ideals, but its only real
prime ideal $0$ is not positive: the involution trace form takes the
values $x,y,-xy$ on $E_{21},E_{31},E_{41}$, and no ordering of $k$ makes
all three positive. Thus its zero ideal is real primitive but is not
the kernel of a $*$-representation over any ordered field extension of
$k$. The strongly anisotropic, nowhere definite form
$\langle1,x,y,-xy\rangle$ is discussed in \cite[p.~68]{LLT}.
\end{remark}

\section{Free associative algebras}\label{sec:free}

This section proves Theorem~\ref{thm:free-main}, except for the
algorithmic assertion in part~(3), proved in
Appendix~\ref{sec:algorithms}. Lewin's theorem splits off the
finite-dimensional semisimple part. For at least two generators,
irreducible completions of partial operators then give separating
finite-dimensional simple quotients.

Fix $F=k\langle x_1,\ldots,x_g\rangle$ with $g<\infty$ over a field $k$.
Throughout this section, Section~\ref{sec:freereal}, and
Appendix~\ref{sec:algorithms}, we identify $F^n$ with the row module
$F^{1\times n}$, with componentwise left multiplication.
Table~\ref{tab:free-results} summarizes the results of this section and
Section~\ref{sec:freereal} for finitely generated submodules of $F^n$;
$n=1$ gives left ideals. Only irreducible points are used in the table.
For points, evaluation maps and their kernels, see
Subsection~\ref{subsec:free-maximal}; for $\rad N$ and
$\sqrt[\mathrm{re}]{N}$, see Subsections~\ref{subsec:prime-radicals}
and~\ref{subsec:realbasic}, respectively.

All assertions in the table concern finitely generated submodules;
they do not assert a property of all submodules of the free module.
In the real column, $k=\mathbb R$ or $\mathbb C$ with its usual involution
and $*$ permutes the letters, except that prime separation allows an
arbitrary involution over any field. The separation and radical results
in the \emph{No involution} column of Table~\ref{tab:free-results} are new.
The extension of real prime separation to arbitrary fields and
involutions is also new. For the real Nullstellensatz and real-radical
bounds, we give direct proofs of results of \cite{CHMN,HKN}.

\begin{table}[htbp]
\centering
\small
\setlength{\tabcolsep}{3pt}
\renewcommand{\arraystretch}{1.25}
\begin{tabular}{@{}>{\raggedright\arraybackslash}p{.20\textwidth}
                    >{\raggedright\arraybackslash}p{.38\textwidth}
                    >{\raggedright\arraybackslash}p{.38\textwidth}@{}}
\hline
 & \textbf{No involution} & \textbf{Real case} \\
\hline
Prime separation
 & F.g. semiprime submodules of $F^n$ are finite intersections of f.g. prime
   submodules (Theorem~\ref{thm:free}).
 & F.g. real submodules of $F^n$ are finite intersections of f.g. real prime
   submodules (Theorem~\ref{thm:freereal}). \\[6pt]
Jacobson property
 & They are also intersections of f.g. maximal submodules (Theorem~\ref{thm:free-jacobson}).
 & They are also intersections of f.g. real maximal submodules (Theorem~\ref{thm:free-real-jacobson}). \\[6pt]
Points
 & Pairs $(\pi,v)$ with $\pi$ a f.d. irreducible representation and
   $0\ne v\in(V_\pi)^n$.
 & Pairs $(\pi,v)$ with $\pi$ a f.d. irreducible $*$-representation and
   $0\ne v\in(V_\pi)^n$. \\[6pt]
Point evaluations
 & \multicolumn{2}{l}{$
   \varepsilon_{\pi,v}(f)=\sum_j\pi(f_j)v_j,\;
   f=(f_1,\ldots,f_n)\in F^n,\;
   v=(v_1,\ldots,v_n)\in(V_\pi)^n.$} \\[6pt]
F.g. maximal submodules
 & Exactly $L_{\pi,v}=\ker\varepsilon_{\pi,v}$ at the points above (Proposition~\ref{prop:free-maximal}).
 & Exactly $L_{\pi,v}=\ker\varepsilon_{\pi,v}$ at the real points above (Corollary~\ref{cor:free-real-maximal}). \\[6pt]
Definition of the radical
 & $\rad N$ is the smallest semiprime submodule of $F^n$ containing $N$.
 & $\sqrt[\mathrm{re}]{N}$ is the smallest real submodule of $F^n$
   containing $N$. \\[6pt]
Nullstellensatz
 & $\rad N$ consists of the vectors vanishing at all points where
   $N$ vanishes (Corollary~\ref{cor:ordinarynss}).
 & $\sqrt[\mathrm{re}]{N}$ consists of the vectors vanishing at all
   real points where $N$ vanishes (Corollary~\ref{cor:realnss}). \\[6pt]
Finite generation of radicals
 & If $N$ is a f.g. submodule of $F^n$, then $\rad N$ is also a f.g.
   submodule of $F^n$; $\rad N/N$ is finite-dimensional (Corollary~\ref{cor:freeradical}).
 & If $N$ is f.g. with generators of degrees $\le d$,
   then $\sqrt[\mathrm{re}]{N}$ is also f.g. with additional generators
   of degrees $<d$ (Theorem~\ref{thm:free-real-radical}). \\
\hline
\end{tabular}
\caption{Separation, evaluations, and radicals over free associative algebras.}
\label{tab:free-results}
\end{table}

\begin{remark}[Reduction to free algebras]\label{rem:presented-quotient-preimage}
Question~\ref{qu:leftdecomposition} reduces to free associative rings:
every unital ring is a quotient of some $\mathbb Z\langle X\rangle$,
and inverse images under a surjection preserve and reflect primality
and semiprimeness of left ideals. Likewise, the question for $k$-algebras
reduces to free $k$-algebras, and the question for finitely generated
$k$-algebras reduces to finitely generated free $k$-algebras. These
reductions require prime separation for all semiprime left ideals of
the free ring or algebra. We prove it only for finitely generated
semiprime left ideals in finitely generated free algebras over fields,
and more generally for finitely generated semiprime submodules of
finite free modules.

A finitely generated left ideal in a finitely presented quotient algebra
need not have a finitely generated inverse image in the presenting free
algebra. Indeed, if $R=F/J$ with $F$ free and $I\subseteq R$ is generated by the
images of $f_1,\ldots,f_r$, its inverse image is
$J+\sum_iFf_i$. Finite presentation of $R$ as an algebra makes $J$
finitely generated as a two-sided ideal, not necessarily as a left
ideal. For example, in $F=k\langle x,y\rangle$ the two-sided ideal
$J=FxF$ is the direct sum $\bigoplus_{m\ge0}Fxy^m$, so it is not
finitely generated as a left ideal, although $F/J\cong k[y]$ is a
finitely presented algebra. Thus finite presentation over $R$ does not
in general give finite presentation over $F$.
\end{remark}

\subsection{Prime separation}

Theorem~\ref{thm:free} answers Question~\ref{qu:leftdecomposition}
affirmatively for finitely generated semiprime left ideals of $F$. More generally,
it expresses every semiprime submodule whose quotient module is finitely
presented as a finite intersection of prime submodules. We first state the
result, then derive it from
Lewin's theorem by splitting off the finite-dimensional part of the module.

\begin{theorem}[Free algebras]\label{thm:free}
Every finitely presented semiprime left $F$-module is a finite direct sum
of prime modules. Consequently, if $N\subseteq M$ and $M/N$ is finitely
presented, then $N$ is semiprime if and only if it is a finite intersection
of prime submodules of $M$.
\end{theorem}

The main technical ingredient of the proof is
Lemma~\ref{lem:freeparts}, which uses Lewin's theorem: every finitely
presented $F$-module contains a free submodule of finite codimension
\cite[Theorem~2]{Lewin}. Note that this free submodule need not be a
direct summand.

\begin{lemma}[Finite-dimensional part]\label{lem:freeparts}
Let $F=k\langle x_1,\ldots,x_g\rangle$, with $g\geq1$.
Every finitely presented semiprime $F$-module $Q$ has a decomposition
$Q=D\oplus K$, where $D$ is finite-dimensional and semisimple, and
$K$ has no nonzero finite-dimensional submodule and, if nonzero, is
faithful and prime.
\end{lemma}

\begin{proof}
By \cite[Theorem~2]{Lewin}, choose a free submodule $U\subseteq Q$
with $Q/U$ finite-dimensional. Let $D$ be the sum of all finite-dimensional
submodules of $Q$. A free module has no nonzero finite-dimensional
submodule: a nonzero vector generates a copy of the infinite-dimensional
domain $F$. Every element of $D\cap U$ lies in $E\cap U$ for some finite sum $E$
of finite-dimensional submodules. Since $E\cap U$ is a finite-dimensional
submodule of $U$, it is zero. Thus $D\cap U=0$, so $D$ embeds into $Q/U$ and is finite-dimensional.

If $D=0$, take $K=Q$. Otherwise put $T=\Ann_F(D)$.
Lemma~\ref{lem:prime-semiprime-annihilators}(2) shows that $T$ is semiprime.
The action map $F\to\operatorname{End}_k(D)$ has kernel $T$, so
$F/T$ embeds into $\operatorname{End}_k(D)$. Thus $F/T$ is
finite-dimensional and semiprime, hence semisimple Artinian.
Semiprimeness of $Q$ gives $D\cap TQ=0$, so the sequence of
$F/T$-modules
\[
 0\longrightarrow D\longrightarrow Q/TQ
 \longrightarrow Q/(D+TQ)\longrightarrow0
\]
is exact. Since $F/T$ is semisimple, it splits. Composing the
resulting retraction $Q/TQ\to D$ with the quotient map gives a
retraction $r:Q\to D$. Thus $Q=D\oplus K$, where $K=\ker r$.
Moreover, $D$ is semisimple as an $F/T$-module.
Every finite-dimensional submodule of $K$ is a finite-dimensional
submodule of $Q$, hence lies in $D$ by its definition, and therefore
in $D\cap K=0$.

For $0\neq L\subseteq K$, the equality $L\cap U=0$ would embed $L$
into $Q/U$, making $L$ finite-dimensional and hence contained in
$D\cap K=0$. Therefore $L\cap U\neq0$, so $K\cap U$ is essential in $K$.
If $K\ne0$, then $K\cap U$ is a faithful prime module: each of its
nonzero submodules contains a copy of $F$, since $F$ is a domain and
$U$ is free. Lemma~\ref{lem:essential-prime} therefore shows that $K$
is faithful and prime.
\end{proof}

\begin{proof}[Proof of Theorem~\ref{thm:free}]
For $g\geq1$, apply Lemma~\ref{lem:freeparts} and decompose the
finite-dimensional semisimple part into simple modules, which are
prime by Lemma~\ref{lem:basic-prime-modules}(1). Together with the nonzero
$K$, they give a finite direct sum of prime modules. Given a finite direct
sum $M/N=\bigoplus_{i=1}^rH_i$ of prime modules,
let $P_i$ be the kernel of $M\twoheadrightarrow M/N\twoheadrightarrow H_i$.
Then $M/P_i\cong H_i$, so $P_i$ is prime, and $\bigcap_iP_i=N$
because all coordinates vanish precisely on $N$. Conversely, intersections of prime
submodules are semiprime by Lemma~\ref{lem:basic-prime-submodules}(2).
For $g=0$, $F=k$ and the assertion is finite-dimensional linear algebra.
The zero module uses the empty-sum and empty-intersection conventions.
\end{proof}

\begin{example}[Finite free modules]\label{ex:finitefree}
Every finitely generated semiprime submodule $N\subseteq F^n$ is a
finite intersection of finitely generated prime submodules, since $F^n/N$
is finitely presented. Indeed, write $Q=F^n/N=\bigoplus_iH_i$ as a
finite direct sum of prime modules, and put $C_i=\bigoplus_{j\ne i}H_j$.
The inverse images $P_i$ of $C_i$ in $F^n$ are prime and have
intersection $N$. Each $C_i$ is finitely generated: projecting a finite
generating set of $Q$ onto $C_i$ gives generators of $C_i$.
If $b_{i1},\ldots,b_{it_i}\in F^n$ lift these generators, then
$P_i=N+\sum_jFb_{ij}$, so $P_i$ is finitely generated.
For $n=1$, this gives the result for finitely generated semiprime left ideals.
More generally, the theorem applies when $M$ is finitely presented and
$N\subseteq M$ is finitely generated: lifts of generators of $N$, together
with a finite set of relations for $M$, give a finite presentation of $M/N$.
\end{example}

\subsection{Maximal submodules and representations}\label{subsec:free-maximal}

To describe points and evaluation, let
$\pi:F\to\operatorname{End}_k(V_\pi)$ be a unital $k$-algebra map with
$0<\dim_kV_\pi<\infty$, and let $v=(v_1,\ldots,v_n)\in(V_\pi)^n$.
Define the evaluation map $\varepsilon_{\pi,v}:F^n\to V_\pi$ and its kernel by
\begin{equation}\label{eq:evaluation}
 \varepsilon_{\pi,v}(f_1,\ldots,f_n)
 =\sum_{j=1}^n\pi(f_j)v_j,
 \qquad L_{\pi,v}=\ker\varepsilon_{\pi,v}.
\end{equation}
A \emph{point} is a pair $(\pi,v)$ with $v\ne0$; it is
\emph{irreducible} if $V_\pi$ has no proper nonzero invariant subspace.
In the real case, a point is \emph{real} when $\pi$ is a
$*$-representation on a positive-definite inner-product space.
A point is a \emph{zero} of $N$ if $\varepsilon_{\pi,v}(N)=0$;
a vector $p\in F^n$ vanishes there if $\varepsilon_{\pi,v}(p)=0$.

We write $V$ for $V_\pi$ when $\pi$ is clear. After choosing a basis,
$\pi$ is the tuple $X_i=\pi(x_i)$, with $\pi(f)=f(X)$.
Irreducibility is over $k$, not necessarily absolute.

\begin{proposition}[Maximal kernels]\label{prop:free-maximal}
For a maximal submodule $L\subset F^n$, the following are equivalent:
\begin{enumerate}
\item $L$ is finitely generated;
\item $L$ has finite codimension over $k$;
\item $L=L_{\pi,v}$ for a finite-dimensional irreducible representation
$\pi$ and a tuple $v\ne0$.
\end{enumerate}
These submodules are classified by such pairs $(\pi,v)$ up to isomorphism:
$L_{\pi,v}=L_{\pi',v'}$ precisely when there is a linear isomorphism $T$
with $T\pi(f)=\pi'(f)T$ for all $f\in F$ and $Tv_j=v'_j$ for every $j$.
\end{proposition}

\begin{proof}
(1)$\Rightarrow$(2). If $L$ is finitely generated, $V=F^n/L$ is
finitely presented and simple. Lewin's theorem gives a free submodule
of finite codimension. For $g\geq1$ it must be zero: otherwise
simplicity makes $V$ free, whereas no nonzero free $F$-module is simple,
since $F$ has a proper nonzero left ideal $Fx_1$. Hence $V$ is
finite-dimensional. For $g=0$ this is immediate.

(2)$\Rightarrow$(3). The action of $F$ on the finite-dimensional simple module
$V=F^n/L$ defines an irreducible representation $\pi$. With $v_j=e_j+L$,
the tuple $v$ is nonzero and $\varepsilon_{\pi,v}(f)=f+L$, so this
is the quotient map and $L=L_{\pi,v}$.

(3)$\Rightarrow$(2). Since $v\ne0$, the map $\varepsilon_{\pi,v}$
is nonzero and hence surjective onto the finite-dimensional simple
module $V$. Thus $L$ has finite codimension (and is maximal).

(2)$\Rightarrow$(1). Every finite-codimension kernel is finitely
generated, even without maximality. To see this directly, choose a basis $s_1,\ldots,s_t$ of
the quotient and lifts $b_1,\ldots,b_t\in F^n$. Write
$e_j+L=\sum_\ell c_{j\ell}s_\ell$ and
$x_i s_\ell=\sum_r a_{i\ell r}s_r$. The finitely many vectors
\[
 e_j-\sum_\ell c_{j\ell}b_\ell,\qquad
 x_i b_\ell-\sum_r a_{i\ell r}b_r
\]
generate $L$. Indeed, modulo their span as an $F$-module, induction
on word length expresses every vector in the scalar span of the
$b_\ell$. This quotient has dimension at most $t$ and maps onto the
$t$-dimensional quotient $F^n/L$, so the map is an isomorphism.
Finally, equal kernels identify both targets with the same quotient
of $F^n$, giving the asserted isomorphism; the converse is immediate.
\end{proof}

For fixed $\pi$, the tuples are identified by $\operatorname{Aut}_F(V)$.
By Schur's lemma this is the group of units of a division algebra;
for algebraically closed $k$ it consists of the nonzero scalars.
If $N$ is generated by $p^{(1)},\ldots,p^{(r)}$, then
$N\subseteq L_{\pi,v}$ exactly when
$\sum_j\pi(p_j^{(\ell)})v_j=0$ for every $\ell$.

\subsection{The Jacobson property}

Separation by finite-dimensional simple modules rests on two lemmas. The
first reduces it to extending finitely many partial operators; it is used
again, with a different extension, in the real case.

For $Q=F^n/N$ and $m\geq0$, let $Q_{\leq m}$ denote the image of the
vectors of degree at most $m$; it is finite-dimensional, and
$x_iQ_{\leq m-1}\subseteq Q_{\leq m}$.

\begin{lemma}[Truncation]\label{lem:truncation}
Let $N\subseteq F^n$ be generated by vectors of degree at most $e\geq1$,
and put $Q=F^n/N$. Let $H$ be a $k$-vector space containing $Q_{\leq e}$,
with $k$-linear operators $X_1,\ldots,X_g$ such that $X_iu=x_iu$ for
every $u\in Q_{\leq e-1}$. Regard $H$ as an $F$-module through the $X_i$.
Then there is an $F$-module map $Q\to H$ restricting to the identity on
$Q_{\leq e}$.
\end{lemma}

\begin{proof}
Define $\tilde\varphi:F^n\to H$ by
$\tilde\varphi(f_1,\ldots,f_n)=\sum_jf_j(X)(e_j+N)$. For a word $w$ with
$|w|\leq e$, induction on $|w|$ gives $\tilde\varphi(we_j)=we_j+N$: if
$w=x_iw'$, then $w'e_j+N\in Q_{\leq e-1}$, so
$\tilde\varphi(we_j)=X_i(w'e_j+N)=x_iw'e_j+N$. Hence $\tilde\varphi$ agrees
with the quotient map on vectors of degree at most $e$. In particular it
kills the generators of $N$, hence $N$, and induces the required map.
\end{proof}

The second lemma supplies irreducible extensions. It requires no involution.

\begin{lemma}[Irreducible completion]\label{lem:free-completion}
Let $g\geq2$, let $V\subseteq W$ be finite-dimensional $k$-spaces, and
let $T_i:V\to W$ be linear maps. Suppose no nonzero subspace
$L\subseteq V$ satisfies $T_iL\subseteq L$ for all $i$.
For every $h\ge\max\{1,2\dim_kW\}$, these maps extend to an irreducible tuple
on an $h$-dimensional $k$-space containing $W$.
\end{lemma}

\begin{proof}
Extend the $T_i$ arbitrarily to $B_i\in\operatorname{End}_k(W)$.
Put $w=\dim_kW$, $m=h-w\ge1$, and $r=m-w\ge0$.
Let $\chi(t)$ be the characteristic polynomial of $B_1$ and set
\[
 f(t)=1+t^r\chi(t),\qquad U=k[t]/(f).
\]
The polynomial $f$ is monic of degree $m$, so $U$ has basis
$1,t,\ldots,t^{m-1}$, with residue classes understood.
On $U$, let $A$ be multiplication by $t$, and define
\[
 J\left(\sum_{j=0}^{m-1}a_jt^j\right)=a_{m-1}\,1.
\]
The pair $(A,J)$ is irreducible. Indeed, a nonzero common invariant
subspace contains a nonzero polynomial $u$ of some degree $d<m$.
Multiplying by $t^{m-1-d}$ gives a polynomial of degree $m-1$ with
nonzero leading coefficient. Applying $J$ therefore shows that the
subspace contains $1$, and applying powers of $A$ gives every basis
vector of $U$.

Since $m\ge w$, choose maps
\[
 R:U\twoheadrightarrow W,\qquad C:W\to U,\qquad\ker C=V.
\]
For $C$, compose the quotient map $W\to W/V$ with any injection
$W/V\hookrightarrow U$. On $H=W\oplus U$, define
\[
 X_1=\begin{pmatrix}B_1&0\\0&A\end{pmatrix},\qquad
 X_2=\begin{pmatrix}B_2&R\\C&J\end{pmatrix},
 \qquad X_i=B_i\oplus0\quad(i\ge3).
\]
These operators extend the $T_i$, since $C$ vanishes on $V$.
The projection $P:H\to U$ along $W$ has the explicit formula
\[
 -X_1^r\chi(X_1)
 =\begin{pmatrix}-B_1^r\chi(B_1)&0\\0&-A^r\chi(A)\end{pmatrix}
 =\begin{pmatrix}0&0\\0&I_U\end{pmatrix}=P.
\]
Here $\chi(B_1)=0$ by Cayley--Hamilton and $f(A)=I_U+A^r\chi(A)=0$.

Let $L\subseteq H$ be invariant under all $X_i$. Since $P$ is a
polynomial in $X_1$, both $P$ and $1-P$ preserve $L$, giving
\[
 L=(L\cap W)\oplus(L\cap U).
\]
The subspace $L\cap U$ is invariant under $A=X_1|_U$ and
$J=PX_2|_U$, hence is either zero or $U$.
In the first case $L\subseteq W$. For $x\in L$, we have
$PX_2x=Cx$, while $X_2x\in L\subseteq W$ implies $PX_2x=0$.
Thus $L\subseteq\ker C=V$. Since $X_i|_V=T_i$, the subspace $L$
is invariant under all $T_i$, so the hypothesis gives $L=0$. In the second case $U\subseteq L$. For $u\in U$, we have
$X_2u=Ru+Ju$, so $(1-P)X_2u=Ru$ because $1-P$ projects onto $W$.
Both $X_2$ and $1-P$ preserve $L$, hence
\[
 W=R(U)=(1-P)X_2(U)\subseteq L,
\]
so $L=H$. This proves irreducibility over $k$.
\end{proof}

\begin{theorem}[Finite-dimensional simple separation]\label{thm:free-jacobson}
Every finitely presented semiprime $F$-module is a subdirect product of
finite-dimensional simple modules. For a finitely generated submodule
$N\subseteq F^n$, the following are equivalent:
\begin{enumerate}
\item $N$ is semiprime;
\item $N$ is a finite intersection of finitely generated prime submodules;
\item $N$ is an intersection of finitely generated maximal submodules;
\item $N$ is an intersection of maximal submodules of finite codimension.
\end{enumerate}
\end{theorem}

\begin{proof}
For $g\geq2$, write a finitely presented semiprime module as $Q=D\oplus K$
by Lemma~\ref{lem:freeparts}. Projections onto the finite-dimensional
simple summands of $D$ detect its nonzero elements, so it remains to
separate elements of $K$. This direct summand is finitely presented
and has no nonzero finite-dimensional submodule. Write $K=F^r/P$
with $P$ finitely generated, and fix $0\ne q\in K$.

Choose $e\ge1$ at least the degrees of a representative of $q$ and
of generators of $P$. Multiplication gives partial operators
$T_i:K_{\le e-1}\to K_{\le e}$. Any common invariant subspace in
$K_{\le e-1}$ is closed under all generators of $F$, hence is a
finite-dimensional submodule of $K$ and must be zero.
Lemma~\ref{lem:free-completion} therefore extends the $T_i$ to an
irreducible tuple on a finite-dimensional space $H\supseteq K_{\le e}$,
which becomes a simple $F$-module. Since the defining relations have
degree at most $e$, Lemma~\ref{lem:truncation} gives an $F$-module map
$K\to H$ fixing $K_{\le e}$. It does not kill $q$ and is surjective,
because its image is a nonzero submodule of the simple module $H$.
These quotients, together with the simple projections from $D$, detect
every nonzero element of $Q$ and give the asserted subdirect product.

For $g=1$, $F=k[x]$, and $D$ is already separated by its simple
summands. The module $K$ is torsion-free, since a nonzero torsion
element would generate a finite-dimensional submodule. Being finitely
generated over the PID $k[x]$, it is therefore free, say $K\cong k[x]^r$.
For $0\ne q=(q_1,\ldots,q_r)$, choose $q_j\ne0$ and an irreducible
polynomial $p$ not dividing $q_j$; this is possible because $k[x]$ has
infinitely many irreducible polynomials. The coordinate map followed
by reduction modulo $p$ gives a surjection $K\to k[x]/(p)$ not killing
$q$, with finite-dimensional simple target. For $g=0$, $F=k$, and
coordinate projections separate every finite-dimensional vector space.

For $Q=F^n/N$, the simple quotients just constructed give
finite-codimension maximal kernels in $F^n$ whose intersection is $N$,
proving (1)$\Rightarrow$(4). Theorem~\ref{thm:free} gives
(1)$\Rightarrow$(2), with finitely generated prime kernels as in
Example~\ref{ex:finitefree}: their quotients are finitely presented
direct summands of $Q$. Proposition~\ref{prop:free-maximal} gives
(3)$\Leftrightarrow$(4). Finally, maximal submodules are prime by
Lemma~\ref{lem:basic-prime-submodules}(1), so each of (2)--(4) implies (1)
by Lemma~\ref{lem:basic-prime-submodules}(2).
\end{proof}

\begin{remark}[Finite and infinite intersections]
Theorem~\ref{thm:free-jacobson} gives a finite intersection of prime
submodules, but the intersection of finitely generated maximal
submodules may have to be infinite. Indeed, a finite such intersection
would embed $F^n/N$ into a finite product of finite-dimensional simple
modules, forcing $F^n/N$ to be finite-dimensional. For example, when
$g\ge1$, the zero left ideal of $F$ is semiprime and finitely generated,
but cannot be a finite intersection of finitely generated maximal
left ideals.
\end{remark}

\subsection{Prime radicals}\label{subsec:prime-radicals}

\begin{corollary}[Finite generation of the prime radical]\label{cor:freeradical}
Let $N$ be a finitely generated submodule of $F^n$, and let $\rad N$ be the
smallest semiprime submodule containing $N$. Then $\rad N/N$ is
finite-dimensional over $k$, $\rad N$ is finitely generated, and
\[
 \rad N=\bigcap\{P\subseteq F^n: P\text{ is a prime submodule and }N\subseteq P\}.
\]
\end{corollary}

\begin{proof}
Intersections of semiprime submodules are semiprime by
Lemma~\ref{lem:basic-prime-submodules}(2), so $\rad N$ exists.

Put $Q=F^n/N$ and $A=\rad N/N$. Since $Q$ is finitely presented,
\cite[Theorem~2]{Lewin} gives a free submodule $U\subseteq Q$ with
$\dim_k(Q/U)<\infty$. If $U\ne0$, it is faithful and prime, since
every nonzero submodule contains a copy of the domain $F$.
By Zorn's lemma,
choose $K\subseteq Q$ maximal with $K\cap U=0$; unions of chains
preserve this condition. The image of $U$ is essential in $Q/K$:
the inverse image of any nonzero submodule strictly contains $K$
and hence meets $U$ nontrivially. This image is isomorphic to $U$,
so Lemma~\ref{lem:essential-prime} makes $Q/K$ faithful and prime.
By Lemma~\ref{lem:basic-prime-submodules}(3), the inverse image of $K$
in $F^n$ is prime, hence semiprime. It contains $N$, so $A\subseteq K$ and
$A\cap U=0$. The latter also holds when $U=0$. Consequently $A$
embeds into the finite-dimensional space $Q/U$.

A finite $k$-basis of $A$ generates it as an $F$-module. Lifting this basis
and adjoining generators of $N$ gives a finite generating set for $\rad N$.
Moreover $Q/A=F^n/\rad N$ is finitely presented and semiprime, so
Theorem~\ref{thm:free} expresses $\rad N$ as the intersection of all prime
submodules containing it. These are exactly the prime submodules containing
$N$, since prime submodules are semiprime by Lemma~\ref{lem:basic-prime-submodules}(1).
\end{proof}

Theorem~\ref{thm:free-main}, together with the
evaluation description in Proposition~\ref{prop:free-maximal}, gives
the following form of Hilbert's Nullstellensatz.

\begin{corollary}[Hilbert's Nullstellensatz]\label{cor:ordinarynss}
For a finitely generated submodule $N\subseteq F^n$,
\[
 \rad N
 =\bigl\{p\in F^n:\varepsilon_{\pi,v}(p)=0\text{ whenever }
 \varepsilon_{\pi,v}(N)=0\bigr\},
\]
where $\pi$ ranges over finite-dimensional irreducible representations
of $F$ over $k$ and $v$ over nonzero tuples in $(V_\pi)^n$.
The field $k$ is arbitrary.
\end{corollary}

\begin{proof}
By Proposition~\ref{prop:free-maximal}, the evaluation kernels in the
statement are maximal, hence semiprime by
Lemma~\ref{lem:basic-prime-submodules}(1). Thus every such kernel
containing $N$ contains $\rad N$. Conversely, $\rad N$ is finitely
generated by Corollary~\ref{cor:freeradical}, so
Theorem~\ref{thm:free-jacobson} and Proposition~\ref{prop:free-maximal}
separate every $p\notin\rad N$ by one of these evaluation kernels.
\end{proof}

Irreducibility is essential here: kernels of arbitrary finite-dimensional
representations need not be semiprime.

\section{Free algebras with involution}\label{sec:freereal}

This section proves Theorem~\ref{thm:free-real-main}, except for the
algorithmic assertion recalled in Appendix~\ref{app:realradicals}.
Part~(1) follows from Theorem~\ref{thm:free-main}(1), since real
modules are semiprime and their submodules are real. Positive
truncations and adjoint extension prove parts~(2) and~(3).

Let $k$ be a field with involution and let
$F=k\langle x_1,\ldots,x_g\rangle$, with $g<\infty$, carry an involution
extending that of $k$. Part~(1) allows any such involution. For
parts~(2) and~(3), take $k=\mathbb R$ with the identity or
$k=\mathbb C$ with conjugation, and suppose $*$ permutes the letters.
Reality cannot in general be replaced by weak reality, even for a
principal left ideal (Example~\ref{ex:freeweak}).

\subsection{Real prime separation}

\begin{theorem}[Free real modules]\label{thm:freereal}
Let $k$ be a field with involution and let
$F=k\langle x_1,\ldots,x_g\rangle$, with $g<\infty$, carry an involution
extending that of $k$. Every finitely presented real left $F$-module is
a finite direct sum of real prime
modules. Consequently, if $N\subseteq M$ is real and $M/N$ is finitely
presented, then $N$ is a finite intersection of real prime submodules of
$M$; if moreover $M$ is finitely generated, these can be chosen finitely
generated.
\end{theorem}

\begin{proof}
A real module is semiprime by Lemma~\ref{lem:apprealbasic}, so
Theorem~\ref{thm:free} writes $Q=M/N$ as a finite direct sum of prime
modules. Each summand is a submodule of $Q$, hence real by
Lemma~\ref{lem:apprealbasic}. The kernels of the coordinate projections are
real prime submodules with zero intersection; take inverse images in $M$.
If $M$ is finitely generated, choose a surjection $F^r\to M$ with $r$
finite. Each prime summand $H_i$ of $Q$ is finitely presented, so the
kernel of the composite surjection $F^r\to M\to H_i$ is finitely
generated. Its image in $M$ is the required prime submodule.
\end{proof}

\subsection{The real Jacobson property}

For the remainder of this section, $k=\mathbb R$ with the identity, or
$k=\mathbb C$ with complex conjugation, and the involution permutes the
letters: $x_i^*=x_{\sigma(i)}$ for an involutive permutation $\sigma$.
We prove that finitely generated real submodules of finite free
modules are intersections of real maximal submodules. The allowed involutions include
symmetric generators ($\sigma=\mathrm{id}$) and generators paired
by $*$. Recall that a finite-dimensional $*$-representation is an
$F$-module structure on a positive-definite inner-product space $V$
with $\pi(f^*)=\pi(f)^*$. The inner product is symmetric bilinear for
$k=\mathbb R$ and Hermitian for $k=\mathbb C$; thus the adjoint identity
is equivalent to $X_i^*=X_{\sigma(i)}$. Such a module is
real: the inner product induces an additive map
$\overline V\otimes_FV\to k$, $u\otimes w\mapsto\langle u,w\rangle$, sending
$\sum_iu_i\otimes u_i$ to $\sum_i\|u_i\|^2$. It is semisimple by
Lemma~\ref{lem:star-decomposition}.

\begin{lemma}[Positive truncation]\label{lem:posttrunc}
Let $N\subseteq F^n$ be a real submodule, $Q=F^n/N$, and $e\geq1$. There is
a positive-definite inner product on $Q_{\leq e}$ such that
\[
 \langle x_iu,w\rangle=\langle u,x_{\sigma(i)}w\rangle
 \qquad(u,w\in Q_{\leq e-1}).
\]
\end{lemma}

\begin{proof}
Let $Z\subseteq\overline Q\otimes_FQ$ be the real span of the elements
$u\otimes w$ with $u,w\in Q_{\leq e}$; it is a finite-dimensional real
vector space. Fix a norm on $Q_{\leq e}$. The set
$\Sigma=\{u\otimes u:\|u\|=1\}$ is compact, so its convex hull is compact.
It does not contain $0$: a relation $\sum_i\lambda_iu_i\otimes u_i=0$ with
$\lambda_i\geq0$ and $\sum\lambda_i=1$ gives
$\sum_i(\sqrt{\lambda_i}u_i)\otimes(\sqrt{\lambda_i}u_i)=0$, so all
$\sqrt{\lambda_i}u_i=0$ by reality, which is impossible. Hence there is a
real linear functional $\varphi$ on $Z$ that is positive on $\Sigma$.
The twist $\tau(u\otimes w)=w\otimes u$ is a well-defined additive
involution of $\overline Q\otimes_FQ$ preserving $Z$ and fixing $\Sigma$;
replace $\varphi$ by $(\varphi+\varphi\circ\tau)/2$.

For $k=\mathbb R$ put $\langle u,w\rangle=\varphi(u\otimes w)$. For
$k=\mathbb C$ put
$\langle u,w\rangle=\varphi(u\otimes w)-i\varphi(u\otimes iw)$. Using
$(\lambda u)\otimes w=u\otimes\bar\lambda w$ and $\varphi\circ\tau=\varphi$,
one checks that this is a Hermitian form, conjugate-linear in the first
variable, with $\langle u,u\rangle=\varphi(u\otimes u)>0$ for $u\neq0$.
Finally, $\overline{x_iu}=\bar u\cdot x_i^*$ in $\overline Q$, so
$\overline{x_iu}\otimes w=\bar u\otimes x_{\sigma(i)}w$. For $u,w\in Q_{\leq e-1}$ both
sides lie in $Z$, which gives the displayed identity.
\end{proof}

\begin{lemma}[Adjoint extension]\label{lem:adjext}
Let $V$ be a finite-dimensional inner-product space, $V'\subseteq V$ a
subspace with orthogonal projection $E$, and $T_i:V'\to V$ linear maps
with $\langle T_iu,w\rangle=\langle u,T_{\sigma(i)}w\rangle$ for
$u,w\in V'$. Then
\[
 X_i=T_iE+(T_{\sigma(i)}E)^*(1-E)
\]
satisfies $X_i|_{V'}=T_i$ and $X_i^*=X_{\sigma(i)}$.
\end{lemma}

\begin{proof}
Set $A_i=T_iE\in\operatorname{End}(V)$. The hypothesis, applied
to $Eu,Ew$, gives $A_i^*E=EA_{\sigma(i)}$.
On $V'$ the second term in $X_i=A_i+A_{\sigma(i)}^*(1-E)$ vanishes,
so $X_i$ extends $T_i$. Moreover
\[
 X_{\sigma(i)}-X_i^*=EA_{\sigma(i)}-A_i^*E=0.
\]
All adjoints here are adjoints of endomorphisms of $V$.
\end{proof}

\begin{theorem}[The real Jacobson property]\label{thm:free-real-jacobson}
Every finitely generated real submodule $N\subseteq F^n$ is an
intersection of finitely generated real maximal submodules,
equivalently of real maximal submodules of finite codimension. Each of
them is the kernel of an evaluation $\varepsilon_{\pi,v}$ at a
finite-dimensional irreducible $*$-representation.
\end{theorem}

\begin{proof}
Fix $q\in F^n\setminus N$, and choose $e\geq1$ bounding the degrees of
$q$ and of finitely many generators of $N$. Put $Q=F^n/N$ and
$V=Q_{\leq e}$, with the positive-definite inner product of
Lemma~\ref{lem:posttrunc}. The partial multiplication maps
$T_i:Q_{\leq e-1}\to V$ satisfy the adjoint identities required by
Lemma~\ref{lem:adjext}, so they extend to operators $X_i$ on $V$ with
$X_i^*=X_{\sigma(i)}$. Since $F$ is free, these operators make $V$ a
finite-dimensional $*$-representation.

With $v_j=e_j+N\in V$, define the $F$-linear evaluation map
\[
 \psi:F^n\longrightarrow V,\qquad
 \psi(f_1,\ldots,f_n)=\sum_j f_j(X)v_j.
\]
By Lemma~\ref{lem:truncation}, $\psi(p)=p+N$ whenever $\deg p\le e$.
Thus $\psi$ kills the generators of $N$, hence $N$, and
$\psi(q)=q+N\ne0$. It is surjective, since every vector of $V$ is the
class of such a $p$.

By Lemma~\ref{lem:star-decomposition}, $V$ is an orthogonal direct sum
of irreducible $*$-representations, and its coordinate projections
are $F$-linear. Choose a summand $S$ on which $\psi(q)$ has nonzero
projection. The composite $F^n\twoheadrightarrow V\twoheadrightarrow S$
has kernel $L$ containing $N$ and excluding $q$. Since $S$ is simple
and real, $L$ is maximal and real; its finite codimension makes it
finitely generated by Proposition~\ref{prop:free-maximal}. As $q$ was
arbitrary, $N$ is the intersection of these kernels.
\end{proof}

The construction needs neither Lewin's theorem nor the completion lemma:
positivity makes the extension problem trivial.

\begin{corollary}[Real maximal kernels]\label{cor:free-real-maximal}
For a real maximal submodule $L\subset F^n$, finite generation and
finite codimension are equivalent, and hold precisely when
$L=L_{\pi,v}$ for a finite-dimensional irreducible $*$-representation
$\pi$ and $v\ne0$, up to the equivalence of
Proposition~\ref{prop:free-maximal}.
\end{corollary}

\begin{proof}
By Proposition~\ref{prop:free-maximal}, finite generation and finite
codimension are equivalent. Under these conditions, apply
Theorem~\ref{thm:free-real-jacobson} to $L$: a separating real
maximal kernel containing $L$ equals $L$. Conversely, such a pair gives a
finite-dimensional simple real quotient.
\end{proof}

\begin{example}[Weak reality and the Jacobson property]\label{ex:freeweak}
Let $F=\mathbb R\langle x,x^*\rangle$ with its usual involution,
$p=x^*x+x-x^*$ and $N=Fp$. Then $N$ is weakly real but not real.
It is an intersection of maximal left ideals of finite codimension,
but not of real maximal left ideals.
\end{example}
\begin{proof}
On $J=Fx\oplus Fx^*$ define the left $F$-linear map
$\phi(ax+bx^*)=a+b(x^*+1)$. Since $\phi(p)=0$, we have
$N\subseteq\ker\phi$; conversely, $a+b(x^*+1)=0$ gives
$ax+bx^*=-bp$, so $\ker\phi=N$.
If $\sum_i f_i^*f_i\in N$, comparison of constant terms shows that
each $f_i$ has zero constant term, so $f_i=a_ix+b_ix^*$. Put
$u_i=\phi(f_i)$. Applying the left $F$-linear map $\phi$ to
$\sum_if_i^*f_i\in N=\ker\phi$ gives
\[
 0=\sum_i f_i^*u_i
   =x^*\sum_i a_i^*u_i+x\sum_i b_i^*u_i.
\]
Comparing first letters shows that both coefficient sums vanish, so
\[
 \sum_i u_i^*u_i
 =\sum_i a_i^*u_i+(x+1)\sum_i b_i^*u_i=0.
\]
The free algebra is formally real, as comparison of highest-degree
coefficients shows. Hence every $u_i=0$, so $f_i\in N$.
On the other hand, $p+p^*=2x^*x$ and $x\notin N$ by degree,
so $N$ is not real. Finally, $N$ is semiprime by
Proposition~\ref{prop:cyclicreal}, so
Theorem~\ref{thm:free-jacobson} gives the asserted intersection
of maximal left ideals; an intersection of real left ideals would be real.
\end{proof}

\subsection{Real radicals}\label{subsec:free-real-radicals}

Let $N\subseteq F^n$ be finitely generated. Its real radical
$\sqrt[\mathrm{re}]{N}$ (Section~\ref{sec:involution}) need not have finite
codimension over $N$ (Example~\ref{rem:radcompare}), but it is finitely generated. Throughout, $g$ counts all letters,
including both letters of each pair interchanged by $*$.

\begin{example}[Comparison with the prime radical]\label{rem:radcompare}
$\rad N/N$ is finite-dimensional, whereas
$\sqrt[\mathrm{re}]{N}/N$ need not be. For example, in
$\mathbb R\langle x,y\rangle$ with symmetric generators, $N=Fx^2$
forces $x\in\sqrt[\mathrm{re}]{N}$, and the classes of $y^jx$
($j\geq0$) are linearly independent modulo $N$.
\end{example}

\begin{theorem}[Finite generation of the real radical]\label{thm:free-real-radical}
Let $N\subseteq F^n$ be generated by vectors of degree at most $d\geq1$.
Then
\[
 \sqrt[\mathrm{re}]{N}=N+F\bigl(\sqrt[\mathrm{re}]{N}\cap F^n_{<d}\bigr),
\]
where $F^n_{<d}$ consists of vectors of degree less than $d$. In
particular $\sqrt[\mathrm{re}]{N}$ is finitely generated, by the given
generators and at most $n\sum_{j=0}^{d-1}g^j$ further vectors of degree
less than $d$.
\end{theorem}

\begin{proof}
Put $M=F^{1\times n}$, $R=\sqrt[\mathrm{re}]{N}$ and $J=N+F(R\cap M_{<d})$,
so $N\subseteq J\subseteq R$ and $J$ is generated in degrees $\le d$.
Suppose $J\ne R$, and choose $p\in R\setminus J$ of least degree $e$;
then $e\ge d$. Put
\[
 V=(M/J)_{\le e},\qquad U=(M/J)_{\le e-1},\qquad W=(M/R)_{\le e},
\]
and let $\alpha\colon V\to W$ be the natural surjection. Minimality of
$e$ means that $J$ and $R$ agree in degrees $<e$; hence
$U\cap\ker\alpha=0$. We will construct a positive-definite
$*$-representation that detects $p+J$, contradicting $p\in R$.

By Lemma~\ref{lem:posttrunc}, applied to the real submodule $R$, $W$
carries a positive-definite inner product with
$\langle x_iu,w\rangle=\langle u,x_{\sigma(i)}w\rangle$ for
$u,w\in(M/R)_{\le e-1}$. Its pullback to $V$ vanishes on $\ker\alpha$.
To make it positive definite, choose a complement $S\supseteq U$ of
$\ker\alpha$ and give $V=S\oplus\ker\alpha$ the orthogonal-sum inner
product making $\alpha|_S$ an isometry and $\ker\alpha$ carrying any
positive-definite inner product. Since $\ker\alpha$ is orthogonal to $U$, the new inner product
preserves all pairings with $U$:
\[
 \langle v,u\rangle_V=\langle\alpha v,\alpha u\rangle_W
 \qquad(v\in V,\ u\in U).
\]
Since $\alpha$ commutes with multiplication, these identities imply the
partial adjoint relations required by Lemma~\ref{lem:adjext} for the
operators $T_i\colon U\to V$, $u\mapsto x_iu$:
\[
 \langle T_iu,w\rangle_V
 =\langle x_i\alpha u,\alpha w\rangle_W
 =\langle\alpha u,x_{\sigma(i)}\alpha w\rangle_W
 =\langle u,T_{\sigma(i)}w\rangle_V\qquad(u,w\in U).
\]
Lemma~\ref{lem:adjext} extends the $T_i$ to operators $X_i$ on $V$ with
$X_i^*=X_{\sigma(i)}$. These give $V$ a new $F$-module structure by
$x_i\cdot v=X_iv$, agreeing with multiplication in $M/J$ on $U$.
Although $V$ need not be a submodule of $M/J$, this new module is real:
the inner product induces an additive map
$\overline V\otimes_F V\to k$, $\bar v\otimes w\mapsto\langle v,w\rangle_V$,
which sends $\sum_i\bar v_i\otimes v_i=0$ to $\sum_i\|v_i\|^2=0$,
forcing every $v_i=0$.

Since all generators of $J$ have degree $\le d\le e$,
Lemma~\ref{lem:truncation} gives an $F$-module map $\psi\colon M/J\to V$
that is the identity on $V$. Let $f\colon M\to M/J\xrightarrow{\psi}V$
be the composite. Since $V$ is
real, $\ker f$ is a real submodule containing $N$, hence contains $R$.
Thus $p\in R$ gives $f(p)=0$, whereas
$f(p)=\psi(p+J)=p+J\ne0$ because $\psi$ is the identity on $V$ and
$p\notin J$. This is the required contradiction.
Thus $J=R$, and a basis of $R\cap M_{<d}$ supplies at most
$\dim_kM_{<d}=n\sum_{j=0}^{d-1}g^j$ additional generators.
\end{proof}

If $N$ is generated by constant vectors ($d=0$), a change of basis in
$k^n$ identifies $F^n/N$ with a finite free module, which is real.
Thus $\sqrt[\mathrm{re}]{N}=N$, and no additional generators are needed.

\begin{remark}[Relation to earlier finite-generation results]\label{rem:earlier-real-radical}
The finite-generation and degree-bound conclusions of
Theorem~\ref{thm:free-real-radical} are known: see
\cite[Corollary~2.6 and Theorem~3.1]{CHMN} for left ideals, and
\cite[Proposition~2.4]{HKN} and \cite[Theorem~7.3]{Nelson} for row
modules. The latter results imply that the real radical is generated
by $N$ and rows of degree $<d$. Our proof obtains this bound directly
from positive truncation and adjoint extension. The algorithm and
its arithmetic assumptions are recalled in Appendix~\ref{app:realradicals}.
\end{remark}

Theorem~\ref{thm:free-real-main}, together with
the evaluation description in Corollary~\ref{cor:free-real-maximal},
gives the following real Nullstellensatz.

\begin{corollary}[Real Nullstellensatz]\label{cor:realnss}
For a finitely generated submodule $N\subseteq F^n$,
\[
 \sqrt[\mathrm{re}]{N}
 =\bigl\{p\in F^n:\varepsilon_{\pi,v}(p)=0\text{ whenever }
 \varepsilon_{\pi,v}(N)=0\bigr\},
\]
where $\pi$ ranges over finite-dimensional (irreducible) $*$-representations
and $v$ over tuples of vectors.
\end{corollary}

\begin{proof}
Each $\ker\varepsilon_{\pi,v}$ is real, so if it contains $N$ it contains
$\sqrt[\mathrm{re}]{N}$. Conversely, $\sqrt[\mathrm{re}]{N}$ is finitely
generated and real, so Theorem~\ref{thm:free-real-jacobson} separates
every $p\notin\sqrt[\mathrm{re}]{N}$ by such an evaluation.
\end{proof}

\begin{remark}[Relation to earlier real Nullstellens\"atze]\label{rem:earlier-real-nullstellensatze}
For $k=\mathbb R$, Helton, Klep and Nelson
\cite[Corollary~1.11]{HKN} prove Corollary~\ref{cor:realnss} for
$k\langle x_1,x_1^\ast,\ldots,x_g,x_g^\ast\rangle$; their
Section~9 treats $k\langle x_1,\ldots,x_g\rangle$ with
$x_i^\ast=x_i$. The first case was proved for left ideals in
\cite[Theorem~1.6]{CHMN} and for row modules in
\cite[Theorem~1.3]{Nelson}.

We now show how an involution permuting the generators reduces to the
first case. An involutive permutation consists of disjoint transpositions and fixed
points. After renaming,
$g=2r+s$ and
\[
 F=k\langle a_1,a_1^\ast,\ldots,a_r,a_r^\ast,y_1,\ldots,y_s\rangle,
 \qquad y_j^\ast=y_j.
\]
Replacing each $y_j$ by two generators $z_j,z_j^\ast$ interchanged
by the involution gives a free $*$-algebra $B$. The $*$-homomorphisms
$\iota:F\to B$ and $\rho:B\to F$ defined by
\[
 \iota(y_j)=(z_j+z_j^\ast)/2,\qquad
 \rho(z_j)=\rho(z_j^\ast)=y_j,
\]
and fixing every $a_i,a_i^\ast$ satisfy $\rho\iota=1_F$.
For finitely generated $N\subseteq F^n$, put $J=B\iota(N)$, applying
the maps componentwise. The submodule
$\rho^{-1}(\sqrt[\mathrm{re}]{N})$ is real by
Lemma~\ref{lem:apprealbasic} and contains $J$. Hence
$p\notin\sqrt[\mathrm{re}]{N}$ implies
$\iota(p)\notin\sqrt[\mathrm{re}]{J}$.
By \cite[Corollary~1.11]{HKN}, choose a finite-dimensional
$*$-representation $\psi$ and tuple $v$ whose evaluation
vanishes on $J$ and detects $\iota(p)$.
Then $\pi=\psi\circ\iota$ gives
$\varepsilon_{\pi,v}=\varepsilon_{\psi,v}\circ\iota$,
which vanishes on $N$ and detects $p$.
Take an irreducible summand if needed.

Our proofs construct separating representations directly on quotient
truncations using Lemmas~\ref{lem:posttrunc} and~\ref{lem:adjext}.
The same lemmas give the radical degree bound, without extending a
positive functional to all of $F$ or using the GNS construction.
Gr\"obner bases enter only in the appendix's algorithms.
In Corollary~\ref{cor:ordinarynss}, irreducibility is essential:
kernels of arbitrary finite-dimensional evaluations need not be semiprime.
\end{remark}

\appendix
\section{Computing radicals over free algebras}\label{sec:algorithms}

This appendix proves the algorithmic assertion in
Theorem~\ref{thm:free-main}(3), making
Corollary~\ref{cor:freeradical} effective. The procedure uses one-sided normal forms
\cite[Section~2.3]{Mora}, \cite{Green}, and radicals of
finite-dimensional algebras \cite[Sections~2.3.1--2.3.2]{Eberly}.
Throughout, $F=k\langle x_1,\ldots,x_g\rangle$, $g\ge1$, and
$N\subseteq F^{1\times n}$ is given by generators of degrees
$\le d$, with $d\ge1$. Put $D=n\sum_{j=0}^{d-1}g^j$.

\subsection{Normal forms and the prime-radical calculation}

Fix a degree-lexicographic order on row monomials. Left reduction of the
input generators gives a finite left Gr\"obner basis of degrees $\le d$:
leading monomials with the same row index are made suffix-incomparable,
so their left multiples cannot cancel. Normal monomials give a basis of
the quotient, and reduction does not increase degree. A nonzero normal
row of degree $\ge d$ generates a free cyclic submodule: multiplying its
leading monomial on the left cannot introduce a forbidden suffix, whose
length is at most $d$. Consequently every finite-dimensional submodule
lies in the span $V$ of the normal monomials of degree $<d$.

\begin{proposition}[The finite-dimensional calculation]\label{prop:radformula}
Let $B$ be the largest finite-dimensional submodule of $Q=F^n/N$.
If $B=0$, then $\rad N=N$. Otherwise let $E$ be its action algebra and
$S=\ker(F\to E/J(E))$. Then
\[
 \rad N/N=B\cap SQ.
\]
\end{proposition}

\begin{proof}
The sum of all finite-dimensional submodules lies in $V$, so it is
finite-dimensional and equals $B$. Put $A=\rad N/N$.
Corollary~\ref{cor:freeradical} gives $A\subseteq B$.
The quotient $Q/SQ$ is a module over the semisimple ring $F/S$,
so it is semiprime. Hence the inverse image of $SQ$ in $F^n$ contains
$\rad N$, giving $A\subseteq SQ$.
Since $B/A$ is a submodule of the semiprime module $Q/A$, its faithful
action algebra is finite-dimensional and semiprime by
Lemma~\ref{lem:prime-semiprime-annihilators}(2), hence semisimple. The
image of $J(E)$ is a nilpotent ideal in this algebra, so it is zero.
Thus $J(E)$ kills $B/A$, and $S$ annihilates it. Semiprimeness of $Q/A$ now gives
$(B/A)\cap S(Q/A)=0$, whence $B\cap SQ=A$.
\end{proof}

\begin{theorem}[Computing the prime radical]\label{thm:radalgorithm}
Over a field with effective arithmetic and an algorithm for radicals of
finite-dimensional algebras, the following procedure computes $\rad N$.
It adjoins at most $D$ generators, all of degree $<d$.
The field hypothesis holds for $\mathbb Q$ and finite fields.
\end{theorem}

\begin{proof}
All module operations use the normal forms just described.
\begin{enumerate}
\item Starting with $V_0=V$, compute
$V_{r+1}=\{v\in V_r:x_iv\in V_r\text{ for all }i\}$ by linear equations.
After at most $D$ strict decreases the stable space is $B$, the largest
submodule contained in $V$. If $B=0$, return the input generators.
\item Compute the algebra $E$ generated by the action matrices on $B$,
and its radical $J(E)$, using \cite[Sections~2.3.1--2.3.2]{Eberly}.
The algebra-generation calculation takes place in
$\operatorname{End}_k(B)$, of dimension at most $D^2$.
\item Compute left generators of $S=\ker(F\to E/J(E))$ from a word basis
$w_1=1,\ldots,w_t$ of $E/J(E)$: if
$x_iw_j=\sum_\ell c_{ij\ell}w_\ell$ there, use
$s_{ij}=x_iw_j-\sum_\ell c_{ij\ell}w_\ell$.
These generate $S$, since they reduce every word to the chosen basis.
\item Reduce a basis $b_1,\ldots,b_s$ of $B$, represented in degree $<d$,
modulo $K=N+SF^n$. The kernel of
$k^s\to F^n/K$, $(a_j)\mapsto\sum_j a_jb_j+K$, gives a basis of
$B\cap SQ$. Adjoin its lifts to the input generators.
\end{enumerate}
Proposition~\ref{prop:radformula} proves correctness. There are at most
$\dim_k B\le D$ added rows, all represented in degree $<d$.
\end{proof}

\subsection{The real-radical algorithm}\label{app:realradicals}

For $k=\mathbb R$, or $k=\mathbb C$ with conjugation, and an involution
permuting the generators, the Real Algorithm of
\cite[Section~3.1]{CHMN}, in the row-module formulation of
\cite[Section~7.5 and Proposition~7.11]{HKN}, computes
$\sqrt[\mathrm{re}]{N}$ by adjoining rows of degree $<d$ from
sum-of-squares witnesses. Every strict update enlarges
$L\cap F^n_{<d}$, so there are at most $D$ strict enlargements,
using the same low-degree space as in
Theorem~\ref{thm:free-real-radical}. With exact semidefinite
feasibility, the procedure is effective for real algebraic coefficients,
or complex coefficients with real algebraic real and imaginary parts;
the feasibility conditions are polynomial equalities and inequalities
over the real algebraic numbers.

If $d=0$, the quotient by $N$ is finite free, so both radicals equal
$N$. If $g=0$, both procedures reduce to linear algebra over $k$.

\section*{Acknowledgements}
In the preparation of this paper the author used the AI assistants Claude (Anthropic) and ChatGPT (OpenAI) for conceptual discussions regarding proofs and literature, and for editing assistance. All arguments and references were independently verified by the author, who takes full responsibility for the content.

\end{document}